%

\magnification=\magstep1 
\input amstex
\documentstyle{amsppt}
\pagewidth{6.5truein}
\pageheight{8.9truein}
\ifx\refstyle\undefinedZQA\else\refstyle{C}\fi

\define\rightblack{ \null\nobreak\hfill$\blacksquare$}
\define\QED{\rightblack\enddemo}
\define\QNED{\rightblack\csname endproclaim\endcsname}
\define\cala{{\Cal A}}
\define\calp{{\Cal P}}
\define\LD{LD}
\define\Lavlaws{LL}
\define\critfunc{{\text{\rm cr}}}
\define\crit#1{{\critfunc({#1})}}
\define\restrict{\restriction}
\define\Lrestrict{\upharpoonleft}
\define\Lequiv#1/{\simeq^{#1}}
\define\tLequiv#1/{\equiv^{#1}}
\define\id{\text{\rm id}}
\define\cf{\operatorname{cf}}
\define\E{{\Cal E}}
\redefine\O{{\Cal O}}
\define\natnum{\omega}
\define\newbmod{\bmod^{\!\!\prime\,\,}}
\define\iffe{\text{iff}}
\define\appless#1#2{#1({<}#2)}
\define\DehFDG{1}
\define\DehSSGL{2}
\define\DehARLDS{3}
\define\DehBGLDO{4}
\define\DouCPAEE{5}
\define\LavLDLFA{6}
\define\LavOAEER{7}
\define\SieFA{8}
\define\WehGP{9}

\topmatter
\title Finite left-distributive algebras and embedding algebras\endtitle
\author Randall Dougherty and Thomas Jech\endauthor
\affil Ohio State University\\Pennsylvania State University\endaffil
\date January 27, 1997 \enddate
\ifx\refstyle\undefinedZQA
   \thanks The first author was supported by NSF grant number
   DMS-9158092 and by a grant from the Sloan foundation.\endgraf
   The second author was supported by NSF grant number
   DMS-8918299\endthanks
\else
   \thanks The first author was supported by NSF grant number
   DMS-9158092 and by a grant from the Sloan foundation.\endthanks
   \thanks The second author was supported by NSF grant number
   DMS-8918299.\endthanks
\fi
\address Department of Mathematics, Ohio State University,
Columbus, OH 43210\endaddress
\email rld\@math.ohio-state.edu \endemail
\address Pennsylvania State University, 215 McAllister Building,
University Park, PA 16802\endaddress
\email jech\@math.psu.edu \endemail
\abstract
We consider algebras with one binary operation $\cdot$ and one generator ({\it
monogenic}) and satisfying the left distributive law $a\cdot (b\cdot
c)=(a\cdot b)\cdot (a\cdot c)$.  One can define a sequence of finite
left-distributive algebras $A_n$, and then take a limit to get an infinite
monogenic left-distributive algebra~$A_\infty$.  Results of Laver and Steel
assuming a strong large cardinal axiom imply that $A_\infty$ is free; it is
open whether the freeness of $A_\infty$ can be proved without the large
cardinal assumption, or even in Peano arithmetic. The main result of this
paper is the equivalence of this problem with the existence of a certain
algebra of increasing functions on natural numbers, called an {\it embedding
algebra}.  Using this and results of the first author, we conclude that the
freeness of $A_\infty$ is unprovable in primitive recursive arithmetic.
\endabstract
\endtopmatter

\document

%

\head
1. Introduction
\endhead

We consider algebras with one binary operation $\cdot$ and one
generator ({\it monogenic}) and satisfying the left distributive law
$a\cdot (b\cdot c)=(a\cdot b)\cdot (a\cdot c)$; in particular, we look for
a representation of the free algebra.

The word problem for the free monogenic left-distributive algebra was
solved by Laver~\cite{\LavLDLFA}
under the assumption of a large cardinal and subsequently by
Dehornoy~\cite{\DehBGLDO} without such an assumption.
Laver's result uses elementary
embeddings from $V_\lambda$ into $V_\lambda$ under the `application' operation
$\cdot$ defined by
$j \cdot k = \bigcup_{\alpha < \lambda} j(k \cap V_\alpha)$. If there
exists such an embedding~$j$ other than the identity, then the algebra $A_j$
generated by $j$ is free.

When the embeddings in $A_j$ are restricted to an initial segment of
$V_\lambda$, they form a finite monogenic left-distributive algebra
\cite{\LavOAEER}, and these finite algebras can be described without reference
to elementary embeddings. In fact, for every $n$ there is a (unique)
left-distributive operation $*_n$ on the set $A'_n = \{1,2,\dots,2^n\}$
such that $a *_n 1 = a+ 1$ for all $a < 2^n$ and $2^n *_n 1 = 1$.

There is a natural way of defining a limit $A_\infty$ of the algebras $A'_n$,
and one can ask whether $A_\infty$ is free. We reduce this problem
to a simple ($\Pi^0_2$) statement of finite combinatorics, and show that
the answer is affirmative provided there exists a nontrivial elementary 
embedding
from $V_\lambda$ into itself. The crucial fact used in the proof is a
theorem of Laver~and~Steel~\cite{\LavOAEER} on critical points of elementary
embeddings.

It is open whether the freeness of $A_\infty$ can be proved without the
large cardinal assumption, or even in Peano arithmetic. The main result of
this paper is the equivalence of this problem with the existence of a
certain algebra of increasing functions on natural numbers.

We introduce {\it embedding algebras}, which are algebras $(A, \cdot)$
of increasing functions $a \colon \natnum \to \natnum$ endowed with a binary
operation $\cdot$.
The axioms for embedding algebras
state that the operation $a\cdot b$
is left distributive and interacts with critical points
(the critical point of a function is the least number moved by the function)
in the expected way.
If a (nontrivial) embedding
algebra~$A$ exists, then $A_\infty$ is free;
conversely, we construct an embedding algebra
under the assumption that $A_\infty$ is free.

The first author proved \cite{\DouCPAEE} that the critical sequence for a
nontrivial elementary embedding~$j$ yields an enumeration of critical
points in $A_j$ that grows faster than any primitive recursive function.
One consequence of the main theorem is that such a fast-growing
function can be defined under the assumption that $A_\infty$ is free.
It follows that the freeness of~$A_\infty$ is unprovable in primitive
recursive arithmetic.

\bigpagebreak
\head
2. The free monogenic left-distributive algebra
\endhead

We consider algebras with one binary operation $\cdot$ generated by a
single generator that we denote by the symbol $1$.  We shall often
write $ab$ instead of $a\cdot b$, and use the convention that
$abc=(ab)c$.

The {\it left distributive law} is the equality
$$
a(bc)=ab(ac).
\tag"(\LD)"
$$
We let $W=W_{\cala}$ be the set of all words built up from $1$ using
the operation $\cdot$, denote by~$\equiv$ (or by~$\equiv_{\cala}$)
the equivalence relation on~$W$ given by
$$
a\equiv b\qquad\iffe\qquad \text{(\LD)}\models a=b,
$$
and let $\cala=W/{\equiv}$ be the free left-distributive algebra on one
generator.

For the rest of this section, let $(A,\cdot)$ be a left-distributive
algebra generated by $1$.  We will summarize the relevant known
results on such algebras.

\definition
{Definition 2.1} We say that $a$ is a {\it left subterm} of $b$, or
$a<_Lb$, if, for some $c_1,\dots,c_k$ ($k>0$),
$b=ac_1\dots c_k$.
\enddefinition

\proclaim
{Lemma 2.2} \roster\item"(i)" If $a<_Lb$ and $b<_Lc$,
then $a<_Lc$.
\item"(ii)" If $a<_Lb$ in $A$, then $ca<_Lcb$ in $A$.
\endroster
\endproclaim

\demo
{Proof} Part (i) is trivial.  For (ii), use distributivity: if
$b=ac_1\dots c_k$, then
\ifnum\the\mag=\magstep1we have \fi
$cb=c(ac_1\dots c_k)=ca(cc_1)\dots(cc_k)$.\QED

\proclaim
{Theorem 2.3 \rm(Dehornoy \cite{\DehSSGL})} For all $a,b\in \cala$, either
$a\equiv b$ or $a<_Lb$ or $b<_La$.
\endproclaim

(This was also proved by Laver~\cite{\LavLDLFA} under the assumption that
$<_L$ is irreflexive.)

The proof of Theorem~2.3 is quite constructive, using several explicit
recursive constructions on words in $W$.
We will outline the proof of this result below.

\proclaim
{Lemma 2.4 \rm(Dehornoy \cite{\DehSSGL})} If the relation $<_L$ on $A$ is
irreflexive, then $(A,\cdot)$
is free and satisfies left cancellation.
\endproclaim

\demo{Proof}
Let $\pi$ be the canonical homomorphism of the free algebra $\cala$ onto $A$.
If $a\ne b$ in $\cala$,
then either $a <_L b$ or $a >_L b$, so either $\pi a <_L \pi b$ or
$\pi a >_L \pi b$, so
$\pi a \ne \pi b$; therefore, $\pi$ is an isomorphism. For left cancellation,
if $a \ne b$ in $A$, then $a<_Lb$ or $b<_La$ by Theorem~2.3, so $ca<_Lcb$ or
$cb<_Lca$, so $ca \ne cb$ by irreflexivity.\QED

\proclaim{Theorem 2.5 \rm(Dehornoy \cite {\DehBGLDO})} There is an algebra
$(A,\cdot)$
on which $<_L$ is irreflexive.
Consequently, the free algebra is linearly ordered
by $<_L$ and satisfies left cancellation.
\QNED

\definition
{Definition 2.6} The {\it depth} of $a\in W$ is
defined recursively as follows:
$$
\gather
\text{depth}(1)=0, \\
\text{depth}(ab)=\max\{\text{depth}(a),\text{depth}(b)\}+1.
\endgather
$$
The {\it herringbone} $u_k$ of depth $k$ is also defined recursively:
$$
\aligned
u_0&=1, \\
u_{k+1}&=1u_k.
\endaligned
$$
\enddefinition

One can also define the
{\it full word} $v_k$, the maximal
word of depth $k$, by $v_0=1$ and $v_{k+1}=v_kv_k.$  Then $v_k$ is
equivalent to $u_k$, because an easy induction shows that
$1v_k=v_{k+1}$.

\proclaim
{Lemma 2.7 \rm(Dehornoy \cite{\DehSSGL, Cor.~2})} If $a$ is a word of
depth $\leq k$, then $au_k=u_{k+1}$ in $\cala$.
\endproclaim

\demo
{Proof} By induction on the depth of $a$ (for all $k$ simultaneously).  For
$a=1$, this is immediate from the definition of $u_k$.  If $a$ has positive
depth, then $a = bc$ where $b$ and $c$ have depth smaller than that of $a$,
and hence ${}\leq k-1$.  Now the induction hypothesis gives
$$
abu_k
=ab(a u_{k-1})
=a(b u_{k-1})
=au_k
=u_{k+1},
$$
as desired.\QED

For $a\in W$, we write $a\rightarrow_{\LD} b$ when $b$ results from $a$ by a
single application of (\LD) from left to right (to a subword of $a$), i.e.,
replacing $x(yz)$ by $xy(xz)$.  We write $a\rightarrow b$ if there is a
sequence $a_0=a,a_1,a_2,\dots,a_k=b$ ($k\ge 0$) such that
$a_i\rightarrow_{\LD}a_{i+1}$ for each $i<k$.

\proclaim
{Proposition 2.8 \rm(Dehornoy \cite{\DehFDG})} There is a mapping
$\partial$ from $W$ to $W$ with the following properties: 
\roster
\item $a \rightarrow \partial a$;
\item if $a\rightarrow_{\LD}b$, then $b\rightarrow \partial a$;
\item if $a\rightarrow b$, then $\partial a \rightarrow \partial b$.
\endroster
\endproclaim

\demo{Proof} First define a binary operation $\otimes$ on $W$ by recursion
on the second argument:
$$\align
a\otimes 1&=a1, \\
a \otimes bc &= (a\otimes b)(a\otimes c).
\endalign$$
(The effect of $a\otimes b$ is to distribute $a$ in $b$ as many times as
possible.)

Then define $\partial$ by another recursion:
$$\align
\partial 1 &= 1,\\
\partial (ab) &= \partial a \otimes \partial b.
\endalign$$
(The word $\partial a$ contains all possible applications of (\LD) within $a$.)

Now everything used here is (or can be viewed as being) defined by recursion,
including~$\rightarrow$ (in terms of $\rightarrow_{\LD}$) and even
$\rightarrow_{\LD}$: $a \rightarrow_{\LD} b$ iff either $a$ has the form
$a_1(a_2a_3)$ and $b = (a_1a_2)(a_1a_3)$, or $a$ and $b$ have the forms
$a_1a_2$ and $b_1b_2$, respectively, and either $a_1 \rightarrow_{\LD} b_1$ and
$a_2 = b_2$, or $a_1 = b_1$ and $a_2 \rightarrow_{\LD} b_2$.  One can now prove
a sequence of statements by straightforward inductions:
$$\allowdisplaybreaks
\setbox0=\hbox{if $b\rightarrow_{\LD} b'$,
   then $a\otimes b \rightarrow_{\LD} a \otimes b'$;}
\edef\tempdimen{\the\wd0}
\define\temp#1{\hbox to\tempdimen{\hss$#1$\hss}}
\alignat 2
\temp{ab \rightarrow a\otimes b;}
&&\qquad&\text{(induct on $b$)} \\
\temp{a\otimes(b\otimes c) \rightarrow (a\otimes b)\otimes(a\otimes c);}
&&\qquad&\text{(induct on $c$)} \\
\temp{\text{if $a\rightarrow a'$, then $a\otimes b \rightarrow a' \otimes b$;}}
&&\qquad&\text{(induct on $b$)} \\
\temp{\text{if $b\rightarrow_{\LD} b'$, then $a\otimes b \rightarrow_{\LD}
   a \otimes b'$;}}
&&\qquad&\text{(induct on $b\rightarrow_{\LD}b'$)} \\
\temp{\text{if $b\rightarrow b'$, then $a\otimes b \rightarrow a \otimes b'$;}}
&&\qquad&\text{(induct on $b\rightarrow b'$)} \\
\temp{a \rightarrow \partial a;}
&&\qquad&\text{(induct on $a$)} \\
\temp{a_1a_2(a_1a_3) \rightarrow a_1 \otimes a_2a_3;}
&&\qquad& \\
\temp{\text{if $a\rightarrow_{\LD} b$, then $b \rightarrow \partial a$;}}
&&\qquad&\text{(induct on $a\rightarrow_{\LD}b$)} \\
\temp{\text{if $a\rightarrow_{\LD} b$, then
   $\partial a \rightarrow \partial b$;}}
&&\qquad&\text{(induct on $a\rightarrow_{\LD}b$)} \\
\temp{\text{if $a\rightarrow b$, then $\partial a \rightarrow \partial b$.}}
&&\qquad&\text{(induct on $a\rightarrow b$)}
\endalignat
$$
This gives the desired properties.\QED

\proclaim{Lemma 2.9 \rm(Dehornoy \cite{\DehSSGL})} If $a<_Lb$ in $W$ (i.e.,
$a$ is a left subterm of $b$ in $W$, with no use of the distributive law), and
$b \rightarrow b'$, then there is a left subterm $a'$ of $b'$ in $W$ such that
$a \rightarrow a'$. \endproclaim

\demo{Proof} A straightforward induction on the length of the derivation
$b \rightarrow b'$. \QED

\demo{Proof of Theorem~2.3} From Proposition~2.8, it follows that, if $a
\equiv_\cala b$, then $a \rightarrow \partial^m b$, whenever $m$ is at least the
length of an (\LD)-derivation of $a \equiv b$.  Now, let $a$ and $b$ be words
in $W$, and choose~$k$ so that both words are of depth $\leq k$.  By
Lemma~2.7, we have $au_k \equiv u_{k+1}$ and $bu_k \equiv u_{k+1}$, so
$au_k\rightarrow \partial^m u_{k+1}$ and $bu_k\rightarrow \partial^m
u_{k+1}$ for some $m$.  By Lemma~2.9, there are left subterms~$a'$ and~$b'$
of $\partial^m u_{k+1}$ such that $a\rightarrow a'$ and $b\rightarrow b'$.
Since $a'$ and $b'$ are left subterms of the same word, we have either $a' =
b'$, $a' <_L b'$, or $b' <_L a'$ in $W$; therefore, either $a \equiv b$, $a
<_L b$, or $b <_L a$ in $\cala$.  \QED

All of the steps in the proof of Theorem~2.3 are accomplished by explicit
recursions and inductions (on terms, (\LD)-derivations, etc.), and it is
easy to see that the recursions are in fact primitive recursions (on the
depths of terms, the lengths of derivations, etc.).  Therefore, Theorem~2.3
can be proved in a very basic theory of arithmetic.  One such theory
is Primitive Recursive Arithmetic (PRA), which is formalized in a
language containing function symbols for all possible function
definitions using the constant $0$, the successor function ${}'$,
composition, and primitive recursion; it has axioms stating that the
function symbols satisfy their definitions, and that $0' \ne 0$, and
a rule of inference allowing induction on quantifier-free formulas.
(See Sieg~\cite{\SieFA} for more details.)  This theory is among
the weakest of the commonly-studied fragments of arithmetic; it
is often referred to as the formal version of what Hilbert meant by
`finitary reasoning.'  It is not hard to show that the methods
used to prove Theorem~2.3 can be formalized in this theory, so
Theorem~2.3 is provable in PRA.

\medpagebreak

Now consider algebras with two binary operations $\cdot$ and $\circ$.
We use the convention $ab\circ c=(ab)\circ c$, $a\circ bc=a\circ(bc)$.
Let $W_{\calp}$ be the set of all words built up from $1$ using both
operations, and let $\calp$ be the free algebra on one generator under
the equivalence
$$
a\equiv_{\calp} b\qquad\iffe\qquad \text{(\Lavlaws)}\models a=b,
$$
where (\Lavlaws) is the following set of axioms (Laver \cite{\LavLDLFA}):
$$
\gathered
a\circ(b\circ c)=(a\circ b)\circ c \\
(a\circ b)c=a(bc) \\
a(b\circ c)=ab\circ ac \\
a\circ b=ab\circ a
\endgathered
\tag"(\Lavlaws)"
$$
Note that (\LD) is a consequence of (\Lavlaws):
$$
a(bc)=(a\circ b)c=(ab\circ a)c=ab(ac).
$$

The motivation for axioms (\Lavlaws) comes from large cardinal theory.
Let $V_\lambda$
be the collection of all sets of rank less than $\lambda,$ where $\lambda$
is a limit ordinal. Under the assumption that there exists a nontrivial
elementary
embedding $j$ from $V_\lambda$ to $V_\lambda$, let us consider the algebra
$(A_j,\cdot)$ generated from $j$ by the operation of {\it application}
$$
j\cdot k = \bigcup_{\alpha < \lambda} j(k \cap V_\alpha)
$$
and the algebra $(P_j,\cdot,\circ)$ generated from $j$ by $\cdot$ and
composition of embeddings. Laver~\cite{\LavLDLFA} shows, among other
things, that
$(A_j,\cdot)$ and $(P_j,\cdot,\circ)$ are respectively the free monogenic
left-distributive algebra and the free monogenic algebra satisfying axioms
(\Lavlaws).

Again, we summarize some known facts about the algebras $(P,\cdot,\circ)$.

Let $P$ be an algebra with one generator $1$ satisfying (\Lavlaws).  Let
$A\subseteq P$ consist of all values in $P$ of words in $W_{\cala}$;
$A$ satisfies (\LD) and is generated by $1$.

Conversely, one can construct an algebra $P$ from an algebra $A$.
The following construction is implicit in Laver~\cite{\LavLDLFA}, and
described explicitly in Dehornoy~\cite{\DehARLDS, Prop.~2}.

\proclaim{Proposition~2.10 \rm(Laver, Dehornoy)} Any algebra
$(A,{\cdot})$ satisfying (\LD) can be extended and expanded to an algebra
$(P,{\cdot},{\circ})$ satisfying (\Lavlaws).\endproclaim

\demo{Proof (sketch)}
Given $(A,{\cdot})$, let $P \supseteq A$ be the set of formal
compositions of one or more elements of $A$, with two such formal
compositions identified if their equality can be deduced from associativity
of composition and the rule $a \circ b = ab \circ a$.  Define $\circ$ and
$\cdot$ for two such compositions $a_1 \circ \dots \circ a_n$ and
$b_1 \circ \dots \circ b_m$ by
$$\gather
(a_1 \circ \dots \circ a_n) \circ (b_1 \circ \dots \circ b_m) =
a_1 \circ \dots \circ a_n \circ b_1 \circ \dots \circ b_m, \\
(a_1 \circ \dots \circ a_n) \cdot (b_1 \circ \dots \circ b_m) =
a_1(\dots(a_n(b_1))\dots) \circ \dots \circ a_1(\dots(a_n(b_m))\dots).
\endgather$$
This is well-defined on $P$ and satisfies (\Lavlaws). \QED

Note that, if $A$ is generated by $1$ using $\cdot$, then $P$ is
generated by $1$ using $\cdot$ and $\circ$.

\proclaim{Lemma 2.11} Every element of the free
(\Lavlaws)-algebra~$\calp$ can be written in the form $a_1 \circ \dots
\circ a_n$ for some $a_1,\dots,a_n \in \cala$. \endproclaim

\demo{Proof}
Induct on the form of $p$ as a word in $W_\calp$.
If $p=1$, we are done.  Otherwise,
$p$ has the form $qr$ or $q \circ r$, where we may assume that
$q = a_1 \circ \dots \circ a_n$ and $r = b_1 \circ \dots \circ b_m$
with $a_i,b_j \in \cala$.  We then have
$$q\circ r = a_1\circ\dots\circ a_n\circ b_1\circ\dots\circ b_m$$ and
$$qr = c_1\circ\dots\circ c_m,$$
where $c_j = a_1(a_2(\dots(a_n(b_j))\dots))$, so $p$ has the desired form.
\QED

In the following proposition, the left-to-right
implication is part of Lemma~3 of Laver~\cite{\LavLDLFA}, while the
right-to-left implication uses Lemma~3.2 of that paper.

\proclaim{Proposition 2.12}
Let\/ $(P,{\cdot},{\circ})$ be an algebra
satisfying (\Lavlaws) and generated by\/~$1$, and let $(A,{\cdot})$ be
the subalgebra of\/ $(P,{\cdot})$ generated by\/~$1$.  Then
$P$ is free (with respect to~(\Lavlaws)) if and only if $A$ is free
(with respect to~(\LD)). \endproclaim

\demo{Proof} First, note that
each term $a \in W_\cala$ is either $1$ or of the (unique)
form $a_1b$ for some~$b$.  The same statement can be made about~$b$, and so
on; we eventually find that each such $a$ has a unique expression of the
form $a_1(a_2(\dots(a_n(1))\dots))$ for some $n \ge 0$ and
$a_1,\dots,a_n \in W_\cala$.

The next fact 
(Laver \cite{\LavLDLFA, Lemma 3.2})
we will use is that, if $n,m \ge 1$, $a_i,b_j \in W_\cala$, and
$$a_1(a_2(\dots(a_n(1))\dots)) \equiv_\cala b_1(b_2(\dots(b_m(1))\dots)),$$
then $$a_1 \circ \dots \circ a_n \equiv_\calp b_1\circ\dots\circ b_m.$$  It
will suffice to show that, if $a_1(a_2(\dots(a_n(1))\dots)) \rightarrow_{\LD}
b_1(b_2(\dots(b_m(1))\dots))$, then $a_1 \circ \dots \circ a_n \equiv_\calp
b_1\circ\dots\circ b_m$, since then one can induct on (\LD)-derivations.
(Note that an application of left distributivity cannot start or finish with
the term~$1$, so no term other than~$1$ is equivalent to~$1$
under~$\equiv_\cala$.) If $a_1(a_2(\dots(a_n(1))\dots)) \rightarrow_{\LD}
b_1(b_2(\dots(b_m(1))\dots))$, then there are two cases:  either the
application of left distributivity occurs within a single term $a_i$, or it
changes $a_i(a_{i+1}(x))$ into $a_ia_{i+1}(a_i(x))$ for some $i$.  In the
first case, we get from $a_1 \circ \dots \circ a_n$ to $b_1\circ\dots\circ
b_m$ by applying left distributivity within $a_i$; in the second case, we get
from $a_1 \circ \dots \circ a_n$ to $b_1\circ\dots\circ b_m$ by replacing $a_i
\circ a_{i+1}$ with $a_ia_{i+1} \circ a_i$.  Both of these changes are
permitted by~(\Lavlaws), so $a_1 \circ \dots \circ a_n \equiv_\calp
b_1\circ\dots\circ b_m$.

We are now ready to show that, if $A$ is free, then $P$ is free.  Assume
$A$ is free, and let $p,q \in W_\calp$ be words such that $p = q$ in $P$;
we must show that $p \equiv_\calp q$.  By Lemma 2.11, there are
$n,m \ge 1$ and $a_i,b_j \in W_\cala$ such that $p \equiv_\calp
a_1 \circ \dots \circ a_n$ and $q \equiv_\calp b_1\circ\dots\circ b_m$.
Since $p=q$ in $P$, $p1 = q1$ in $P$, so
$(a_1 \circ \dots \circ a_n)\cdot1 = (b_1\circ\dots\circ b_m)\cdot1$ in $P$,
so $a_1(a_2(\dots(a_n(1))\dots)) =
b_1(b_2(\dots(b_m(1))\dots))$ in $P$ and hence in $A$.  Since
$A$ is free, we have $a_1(a_2(\dots(a_n(1))\dots)) \equiv_\cala
b_1(b_2(\dots(b_m(1))\dots))$.  Now the preceding paragraph gives
$a_1 \circ \dots \circ a_n \equiv_\calp
b_1\circ\dots\circ b_m$, so $p \equiv_\calp q$, as desired.

Now assume that $P$ is free; we must show that $A$ is free.  To do this,
we will show that, if $a,b \in W_\cala$ and $a \not\equiv_\cala b$, then
$a\ne b$ in $A$.  By Proposition~2.10, there is an algebra $P'$ extending
the free algebra $\cala$ which satisfies~(\Lavlaws).  Since
$a \not\equiv_\cala b$, we have $a \ne b$ in $P'$, so
$a \not\equiv_\calp b$.  Since $P$ is free, $a \ne b$ in $P$ and hence
in $A$.  Therefore, $A$ is free. \QED

It is not hard to see that the proof of Propostion~2.12 can be carried
out in PRA; one merely has to use the proof of Proposition~2.10 rather
than the proposition itself when showing ``if $a \not\equiv_\cala b$,
then $a \not\equiv_\calp b$.''

Now consider the algebras $A_j$ and $P_j$ of elementary embeddings. For each
nontrivial elementary embedding from $V_\lambda$ to itself, let $\crit a$
be the {\it critical point} of $a$, the least ordinal moved by $a$. Let
$\Gamma$ be the set of all critical points of elements of $A_j$.
We note that
$$
\crit{ab} = a(\crit b),\qquad \crit{a \circ b} = \min(\crit a,\crit b).
$$
Consequently, the critical point of every $a \in P_j$ is in $\Gamma$, and
every $a \in P_j$ maps $\Gamma$ into $\Gamma.$

\proclaim {Theorem 2.13 \rm(Laver and Steel \cite{\LavOAEER})}
The set $\Gamma$ has order type $\omega$. \QNED

\proclaim {Theorem 2.14 \rm(Laver \cite{\LavOAEER})}
For every $a,b \in A_j$, if $a \ne b$, then $a(\gamma)\ne
b(\gamma)$ for some $\gamma \in\Gamma$. \QNED

Let $\kappa_0$ be the critical point of $j$, and, for all $n$, let $\kappa_{
n+1} = j(\kappa_n)$.

\proclaim
{Lemma 2.15}
\roster
\item"(i)" If $a\in A_j$ has depth at most $n$, then $a(\kappa_n)=\kappa_{
n+1}.$

\item"(ii)" For every $a\in P_j$, there are natural numbers $d >0$ and $N$
such that
$a(\kappa_n)=\kappa_{n+d}$ for all $n \geq N$.
\endroster
\endproclaim

\demo{Proof} (i) By induction on the depth of $a$:
$$
ab(\kappa_n)=ab(a(\kappa_{n-1}))=a(b(\kappa_{n-1}))=a(\kappa_n)=\kappa_{n+1}.
$$

(ii) By Lemma 2.11, we have $a=a_1\circ\dots\circ a_d$ for some
$a_1,\dots,a_d\in A_j$.
\QED

To conclude this section, we remark that one can adjoin to $P_j$ the identity
embedding~$\id$. The extended algebra still satisfies axioms (\Lavlaws), as
well as these rules:

$$
\id\cdot a = a,\quad a\cdot\id=\id,\quad a\circ\id=
\id\circ a=a.
$$

\bigpagebreak
\head 3. A sequence of finite algebras \endhead

In this section, we will construct, for each natural number~$n$, an
algebra $A'_n$ on
\ifnum\the\mag=\magstep1the set \fi
$\{1,2,\dots,2^n\}$ with a binary operation
$*_n$ satisfying the left distributive law.  We will then construct a
second operation $\circ_n$ on this set so that the resulting
two-operation algebra $P'_n$ satisfies~(\Lavlaws).  The subscripts on the
operations will sometimes be omitted while a fixed~$n$ is being
considered.

The construction of these algebras is due to Laver; Wehrung proved some
additional properties of them.  The proof of
the following theorem has been reconstructed independently by
several people, including the authors; the presentation
here is similar to that of Wehrung~\cite{\WehGP}.
(See also Dehornoy \cite{\DehARLDS,~Prop.~7}.)

\proclaim{Theorem 3.1${}'$ \rm(mostly Laver)}
Let $n \ge 0$.

\roster
\item"(a)" 
There is a unique
left-distributive operation $*_n$ on $\{1,2,\dots,2^n\}$ such that
$$
a *_n 1= a+1 \text{ for all } a<2^n, \text{ and } 2^n *_n 1=1.
$$

\item"(b)" 
There is a unique
additional operation $\circ_n$ on $\{1,2,\dots,2^n\}$
such that $*_n$ and $\circ_n$ satisfy axioms (\Lavlaws).
\endroster
\endproclaim

The operation $*_n$ is defined by double recursion; $a *_n b$ is defined
by an outer descending recursion on~$a$ and an inner ascending recursion
on~$b$.  The recursive formulas are as follows:
$$ 2^n *_n b = b;\tag3.1a$$
if $a<2^n$, then
$$a *_n 1 = a+1;\tag3.1b$$
if $a<2^n$ and $b<2^n$, then
$$a *_n (b+1) = (a*_nb)*_n(a+1).\tag3.1c$$
In order to see that this is a valid recursion, we must maintain
the inductive condition
$$a*_nb > a \qquad\text{if $a<2^n$.}\tag3.2${}'$ $$
This clearly holds for $a*_n1$.  For $a*_n(b+1)$ with $a<2^n$, we have
$a*_nb > a$ by the induction hypothesis, so $(a*_nb)*_n(a+1)$ has
already been defined.  If $a*_nb = 2^n$, then $a*_n(b+1) = 2^n*_n(a+1) =
a+1 > a$; if $a*_nb < 2^n$, then $a*_n(b+1) = (a*_nb)*_n(a+1) > a*_nb >
a$.  Therefore, (3.2${}'$) holds for $a*_n(b+1)$ as well, so the recursion
can continue.

The equations~(3.1) can be deduced from left distributivity and the
equations $a *_n 1= a+1$ ($a<2^n$) and $2^n *_n 1=1$.  This is obvious
for (3.1b); for (3.1c) and (3.1a), we have
$$\gather
a *_n (b+1) = a *_n (b *_n 1) = (a*_nb)*_n(a *_n 1) = (a*_nb)*_n(a+1), \\
2^n *_n b = 2^n *_n (\underbrace{1 *_n \dots *_n 1}_{\text{$b$ times}}) =
\underbrace{(2^n *_n 1) *_n \dots *_n (2^n *_n 1)}_{\text{$b$ times}} =
\underbrace{1 *_n \dots *_n 1}_{\text{$b$ times}} = b.
\endgather$$
This proves the uniqueness part of Theorem~3.1${}'$(a).

An easy induction on $b$ shows that the equations~(3.1) hold even
when $a = 2^n$, if we treat addition as being modulo~$2^n$.  (Since we
are working with the set $\{1,2,\dots,2^n\}$, it will be convenient to
treat reduction modulo $2^n$ as a mapping into this set;
we will write ``$x \newbmod 2^n$'' to mean the unique
member of $\{1,2,\dots,2^n\}$ which is congruent to~$x$ modulo~$2^n$.
In particular, $0 \newbmod 2^n$ will be~$2^n$.)  We will soon show that
the equations also hold for $b = 2^n$, and prove several other
useful properties of~$A'_n$ at the same time.

For any fixed $a$, consider the sequence $a *_n1, a*_n2,\dots,a*_n2^n$
in~$A'_n$.  If $a=2^n$, this sequence is just $1,2,\dots,2^n$.  If
$a<2^n$, then the sequence begins with~$a+1$, and
(by~(3.1c)) each member is
obtained from its predecessor by operating on the right by~$a+1$; hence,
by~(3.2${}'$), the sequence must be strictly increasing as long as its
members remain below~$2^n$.  Once $2^n$ is reached (as must happen in at
most $2^n-a$ steps), the next member will be $a+1$ again, and the
sequence repeats.  Therefore, the sequence $a *_n1, a*_n2,\dots,a*_n2^n$
is periodic (as long as it lasts); each period is strictly increasing
from~$a+1$ to~$2^n$.  We will refer to the number of terms in each
period of this sequence as {\it the period of~$a$ in~$A'_n$}.  (The
period of~$2^n$ in~$A'_n$ is~$2^n$.)

\proclaim{Proposition 3.2}

(a) The period of any~$a$ in~$A'_n$ is a power
of~$2$; equivalently, $a*_n2^n = 2^n$ for all~$a$.

(b) The formulas (3.1) hold modulo\/~$2^n$ in $A'_n$
even when $a$ or~$b$ is $2^n$.

(c) Reduction modulo $2^n$ is a homomorphism from $A'_{n+1}$ to~$A'_n$:
$$(a *_{n+1} b) \newbmod 2^n = (a
\newbmod 2^n) *_n (b \newbmod 2^n)$$ for all $a,b$ in $A'_{n+1}$.

(d) For any $a<2^n$ in $A'_n$, if $p$ is the period of $a$ in $A'_n$, then the
period of~$a+2^n$ in~$A'_{n+1}$ is also $p$, and the period of~$a$
in~$A'_{n+1}$ is either $p$ or $2p$. The period of $2^n$ in $A'_{n+1}$ is $2^n$.
\endproclaim

\demo{Proof} By simultaneous induction on $n$.  Part (a) for $n=0$ is
trivial.

Suppose (a) holds for~$n$.  We noted before that the formulas (3.1)
hold modulo~$2^n$ when $a=2^n$.  If $a<2^n$ but $b=2^n$, then $(b+1)
\newbmod 2^n = 1$ and $a*_n1 = a+1$, while $a*_nb = 2^n$ by (a), and
$2^n*_n(a+1) = a+1$, so~(3.1c) holds even in this
case.  Therefore, (b) holds for $n$.

Part (c) for $n$
is proved by induction, downward on $a$ and upward on $b$, as in the
definition of~$*_{n+1}$.  If $a=2^{n+1}$, then both sides are equal to
$b \newbmod 2^n$.  If $b=1$, then both sides are equal to $(a+1)
\newbmod 2^n$.  If $a<2^{n+1}$ and $b>1$, then the left side is equal to
$$((a \newbmod 2^n) *_n ((b-1) \newbmod 2^n)) *_n ((a+1) \newbmod 2^n)$$
by the induction hypothesis, and the right side is also equal to this
value by~(b).  Therefore, (c) holds for~$n$.

Next, consider (d).  Clearly the period of $2^{n+1}$ in $A'_{n+1}$ is $2^{n+1}$,
twice the period of~$2^n$ in~$A'_n$.  Now suppose $a<2^n$, and let $p$ be
the period of $a$ in $A'_n$.  By (c), for each~$b$ in~$A'_n$, $a*_{n+1}b$
and $(a+2^n)*_{n+1}b$ must each be equal to either $a*_nb$ or
$(a*_nb)+2^n$; if $a*_nb < 2^n$, then both of these values are less than
$2^{n+1}$.  It follows that the periods of $a$ and $a+2^n$ in~$A'_{n+1}$
are at least $p$.  Furthermore, by~(3.2${}'$), we must have
\nopagebreak $(a+2^n)*_{n+1}b
> a+2^n$, so $(a+2^n)*_{n+1}b$ must be equal to $(a*_nb)+2^n$ for all
such $b$, so, in particular, $(a+2^n)*_{n+1}p = 2^{n+1}$; hence, the
period of $a+2^n$ in $A'_{n+1}$ is exactly $p$.  (The same argument
shows that the period of $2^n$ in $A'_{n+1}$ is $2^n$.)  For the period
of $a$ in $A'_{n+1}$, there are two cases.  If $a*_{n+1} p = 2^{n+1}$,
then the period of $a$ in $A'_{n+1}$ is $p$, and we are done.  If not,
$a*_{n+1}p$ must be $2^n$.  Then $a*_{n+1}(p+1)$ must be either $a+1$ or
$a+1+2^n$ by (c), and it must be greater than $2^n$ because $a*_{n+1}b$
increases with $b$ until it reaches $2^{n+1}$, so we must have
$a*_{n+1}(p+1) = a+1+2^n = (a*_n1)+2^n$.  Similarly, using part (c)
along with (3.1c) and (3.2${}'$), we see that $a*_{n+1}(p+b) = (a*_nb)+2^n$
successively for $b=2,3,\dots,p$.  In particular, $a*_{n+1}b < 2^{n+1}$
for $b<2p$ and $a*_{n+1}2p = 2^{n+1}$, so the period of $a$ in $A'_{n+1}$
is $2p$.  This completes the proof of (d) for $n$.

Finally, (a) for $n+1$ (in the first phrasing) follows immediately from
(a) and (d) for $n$.  This completes the induction. \QED

Given these properties of $A'_n$, the proof that the left distributive
law holds in $A'_n$ is a straightforward triple induction (downward on
$a$ and $b$, upward on $c$):
$$2^n*(b* c)=b*c = (2^n* b)*(2^n* c);$$
$$a*(2^n* c)=a* c=2^n*(a* c)=(a* 2^n)*(a* c);$$
if $a,b<2^n$, then
$$a*(b* 1)=a*(b+1) =(a* b)*(a+ 1) = (a*b)*(a*1);$$
and, furthermore, if $c<2^n$, then
$$ \alignat 2
a*(b*(c+1))&=a*((b* c)*(b+ 1)) &&\\
&=(a*(b* c))*(a*(b+1))&&\qquad [b* c>b] \\
&=((a* b)*(a* c))*((a* b)*(a+ 1)) &&\\
&=(a* b)*((a* c)*(a+1))&&\qquad [a*b>a] \\
&=(a* b)*(a*(c+1)).&&
\endalignat $$

We now want to define a second operation $\circ = \circ_n$ so that the
resulting algebra $$P'_n = (\{1,2,\dots,2^n\},*_n,\circ_n)$$ satisfies
Laver's axioms (\Lavlaws).  In particular, it will have to be true that
$(a\circ_n b)*_n 1 = a *_n (b*_n 1)$; therefore, we must define
$$a\circ_n b=(a*_n(b+1))-1,$$
where the addition and subtraction
are performed modulo $2^n$.  (So we immediately get the uniqueness in
Theorem~3.1${}'$(b).) This definition makes it immediate that
reduction modulo $2^n$ is a homomorphism from $P'_{n+1}$ to $P'_n$.
We now proceed to prove the four laws (\Lavlaws).  All addition and
subtraction below is modulo $2^n$.

First, one can show that $2^n\circ x = x \circ 2^n = x$ as follows:
$$\gather
2^n \circ x = 2^n * (x+1) - 1 = (x+1) - 1 = x, \\
x \circ 2^n = x * (2^n+1) - 1 = (x*1) - 1 = x.
\endgather$$

The proof of $(a\circ b)*c=a*(b*c)$ is by induction on $c$:
$$(a\circ b)* 1=(a\circ b)+1=a*(b+1)=a*(b* 1);$$
$$\multline(a\circ b)*(c+1)=((a\circ b)*c)*((a\circ b)+1)=(a*(b*c))*(a*(b+1))
\\=a*((b*c)*(b+1))=a*(b*(c+1)).\endmultline$$

Next, $a\circ b=(a*b)\circ a$ because
$$(a\circ b)+1=a*(b+1)=(a*b)*(a+1) = ((a*b)\circ a)+1.$$

The proof of the associative law $a\circ(b\circ c)=(a\circ b)\circ c$
is as follows:
$$\aligned
(a\circ(b\circ c))+1&=a*((b\circ c)+1) \\
&=a*(b*(c+1)) \\
&=a*((b*c)*(b+1)) \\
&=(a*(b*c))*(a*(b+1)) \\
&=((a\circ b)*c)*((a\circ b)+1)\\
&=(a\circ b)*(c+1) \\
&=((a\circ b)\circ c)+1.\\
\endaligned
$$

Finally, to prove that $a*(b\circ c)=(a*b)\circ (a*c)$,
proceed by induction downward on $b$.
For $b=2^n$, we have
$$a*(2^n\circ c)=a*c = 2^n\circ(a*c)=(a*2^n)\circ(a*c).$$
If $b<2^n$, then
$$
\align
a*(b\circ c)&=a*((b*c)\circ b) \\
&=(a*(b*c))\circ (a*b) \qquad\qquad [b*c>b] \\
&=((a*b)*(a*c))\circ (a*b) \\
&=(a*b)\circ (a*c).
\endalign
$$
This completes the proof that $P'_n$ satisfies (\Lavlaws), so Theorem~3.1${}'$
is proved.

\medpagebreak

The following fact will be useful later:
$$\text{if}\quad a\ne 2^n \quad\text{or}\quad b \ne 2^n, \quad\text{then}\quad
a\circ b \ne 2^n.\tag3.3${}'$ $$
This is proved by cases.  If $a \ne 2^n$, then $a*(b+1)>a$ by (3.2${}'$),
so $a*(b+1) \ne 1$, so $a \circ b \ne 2^n$.  If $a=2^n$ but $b \ne 2^n$,
then $a\circ b = b \ne 2^n$.

We remark that Theorem~3.1${}'$ can be rephrased
slightly, replacing $2^n$ by $0$:

\proclaim{Theorem 3.1 \rm(same credits as for 3.1${}'$)}
There are unique operations $*_n$ and $\circ_n$ on $A_n=P_n=\{0,1,\dots,2^n-1\}$
such that the axioms (\Lavlaws) hold and, for all $a \in P_n$,
$$
a*_n 1 = a+1 \bmod 2^n.
$$
\QNED

This has no effect on the structure of the algebras, but it affects statements
referring to the ordering of the elements of the algebra.  In particular,
(3.2${}'$) and (3.3${}'$) become:
$$\gather
\text{either}\quad a*_nb = 0 \quad \text{or} \quad a*_nb > a;\tag3.2\\
\text{if}\quad a\ne 0 \quad\text{or}\quad b \ne 0, \quad\text{then}\quad
a\circ b \ne 0.\tag3.3
\endgather$$
Also, the ordinary $\bmod$ operation now gives the homomorphism from
$P_{n+1}$ to $P_n$.

The element $0$ (or $2^n$) of the algebra plays the role that the
identity embedding played at the end of section~2: 
$$
0*a=a,\quad a*0=0,\quad a\circ 0=0\circ a=a.
$$

\bigpagebreak
\head
4. The limit algebras $A_{\infty}$ and $P_{\infty}$
\endhead

Using the finite algebras $(A_n,*_n)$ and $(P_n,*_n,\circ_n)$, we
construct monogenic algebras $(A_{\infty},\cdot)$ and
$(P_{\infty},\cdot,\circ)$.  Let $W_{\cala}\subset W_{\calp}$ be the
sets of words built up from $1$ using $\cdot$ and using~$\cdot,\circ$,
respectively.  The set of all positive integers can be embedded in
$W_{\cala}$ by identifying each positive integer with a word in $W_{\cala}$, by
recursion:
$$
\gather
1=1, \\
a+1=a\cdot 1.
\endgather
$$
We also adjoin $0$ to $W_{\calp}$, letting
$W_{\calp}^*=W_{\calp}\cup\{0\}$ and $W_{\cala}^*=W_{\cala}\cup\{0\}$,
and add rules
$$
0\cdot a=a,\quad a\cdot0=0,\quad a\circ 0=0\circ a=a.
$$
For every word $a\in W_{\calp}^*$ and every $n\geq 0$, let
$[a]_n$ be the value of $a$ in $P_n=\{0,1,\dots,2^n-1\}$,
and consider the equivalence relation $\equiv_\infty$ defined by:
$$
a \equiv_\infty b \quad \iffe \quad [a]_n=[b]_n\text{ for all }n\geq 0.
$$
We let $A_{\infty}$ and $P_{\infty}$ be, respectively, the
quotients by $\equiv_\infty$ of $W_{\cala}$ and $W_{\calp}$.
Clearly $A_\infty$ and~$P_\infty$ are generated by $1$; also, they
satisfy (\LD) and (\Lavlaws), respectively, because $A_n$ and~$P_n$ do.
(In fact, an equivalent definition for $A_\infty$ and $P_\infty$ is that
they are the subalgebras generated by $1$ of the inverse limits of
the algebras $A_n$ and~$P_n$, respectively.)
Moreover, $A_\infty\subseteq P_\infty$.  We shall investigate the possibility
that $A_\infty$ or $P_\infty$ is free.

\proclaim
{Lemma 4.1} For every $a\in W_P$ and every $n$,
$[a]_{n+1}$ is either $[a]_n$ or $[a]_n+2^n$.
\endproclaim

\demo
{Proof}  This follows immediately from the fact that reduction modulo $2^n$
is a homomorphism from $P_{n+1}$ to $P_n$.  \QED

Note that, as a corollary, if $[a]_n\neq 0$, then $[a]_{n+1}\neq 0$.

\definition
{Definition 4.2} Let $a\in W_{\calp}$ be such that
$[a]_n\neq 0$ for some $n$.  The {\it signature}
$s(a)$ of $a$  is the largest $n$ such that $[a]_n=0$.
\enddefinition

By Lemma 4.1, 
for each $n > s(a)$, $2^{s(a)}$ is
the largest power of $2$ which divides $[a]_n$.

\proclaim
{Lemma 4.3} Let $a,b\in W_{\calp}$ be such that
$[b]_n\neq 0$ for some $n$.  Then, for every $n\geq 0$,
$$
[ab]_n=0\qquad\iffe\qquad [a\cdot 2^{s(b)}]_n=0.
$$
\endproclaim

\demo
{Proof} If $[a\cdot 2^{s(b)}]_n=0$, then $[a]_n *_n [2^{s(b)}]_n = 0$, so
$2^{s(b)}$ is a multiple of the period of $[a]_n$ in $P_n$.  But
$[b]_n$ is a multiple of $2^{s(b)}$, so $[a]_n*_n [b]_n=0$, so
$[ab]_n=0$.

On the other hand, suppose $[ab]_n=0$; then $[a]_n*_n [b]_n=0$.
If $s(b) \ge n$, then $[2^{s(b)}]_n=0$, so $[a]_n*_n[2^{s(b)}]_n=0$.
If $s(b) < n$, 
let $q$ be the period of $[a]_n$ in $A_n$;  then
$q$ divides $[b]_n$, and since $q$ is a power of 2, $q$ divides
the largest power of $2$ dividing $[b]_n$, which is
$2^{s(b)}$.  This again gives $[a]_n*_n[2^{s(b)}]_n=0$. Hence, in either
case, $[a\cdot 2^{s(b)}]_n=0$.  \QED

\proclaim
{Corollary}  $s(ab)=s(a\cdot 2^{s(b)})$.
\endproclaim

\proclaim
{Theorem 4.4}  The following are equivalent:
\roster
\item"(i)"  $(A_{\infty},\cdot)$ is free.
\item"(ii)"  $(P_{\infty},\cdot,\circ)$ is free.
\item"(iii)"  $A_{\infty}$ satisfies the left cancellation law.
\item"(iv)"  $<_L$ on $A_{\infty}$ is irreflexive.
\item"(v)"  If $a<_Lb$ in $A_{\infty}$, then $[a]_n<[b]_n$ for
all but finitely many $n$.
\item"(vi)"  For every $a\in W_{\cala}$,
there is an $n$ such that $[a]_n\neq 0$.
\item"(vii)"  For every $k\geq 0$, there is an $n$ such that $[u_k]_n\neq 0$.
\item"(viii)"  For every $k\geq 1$, there is an $n$ such that
$[1\cdot k]_n\neq 0$.
\endroster
\endproclaim

\demo
{Proof}
\roster
\item"(i)$\leftrightarrow$(ii):" Proposition 2.12.
\item"(i)$\rightarrow$(viii):"  Assume that, for all $n$, $[1\cdot k]_n=0$.
Then, in each $A_n$, $1*(k+1)=
(1* k)* 2=0* 2=2=1* 1$.  However, it is easy to see that the
word $1\cdot 1$ is inequivalent in the free algebra
to any other word, because no application of the distributive law
can start from or result in $1\cdot1$. Therefore, $A$ is not free.
\item"(viii)$\rightarrow$(vii):"  By induction on $k\geq 0$, we prove that
$[u_k]_n\neq 0$ for some $n$.  Assume that this is true for $k$, and let
$s=s(u_k)$ be the signature of $u_k$.  Let $n$ be such that
$[1\cdot 2^s]_n\neq 0$.
By Lemma 4.3, we have $[u_{k+1}]_n=[1\cdot u_k]_n\neq 0$.
\item"(vii)$\rightarrow$(vi):"  Let $k$ be the depth of $a$.
We show that, if $[a]_n=0$, then $[u]_n=0$, where $u=u_k$.
By Lemma~2.7, $au=u_{k+1}=uu$.  If $[a]_n=0$, then $[au]_n=0* [u]_n=[u]_n$;
since $[u]_n* [u]_n$ is either $0$ or ${>}[u]_n$ by formula (3.2), we
have $[u]_n=0$.
\item"(vi)$\rightarrow$(v):"  If $[a]_m\neq 0$ for some $m$,
then $[a]_n\neq 0$ for all $n\ge m$.  Suppose $a<_Lb$, say $b=ac_1\dots c_k$.
Let $n$ be sufficiently large that $[a]_n\neq 0$,
$[ac_1]_n\neq 0$, $[ac_1 c_2]_n\neq 0,\dots,[b]_n\neq 0$.
By (3.2), we have $[a]_n<[ac_1]_n<\dots<[b]_n$.
\item"(v)$\rightarrow$(iv):"  Trivial.
\item"(iv)$\rightarrow$(iii):"  Lemma 2.4.
\item"(iii)$\rightarrow$(viii):"  As for (i)$\rightarrow$(viii),
if $[1\cdot k]_n=0$ for all $n$,
then $1\cdot(k+1)=1\cdot 1$ in $A_\infty$, violating left cancellation.
\item"(iv)$\rightarrow$(i):"  Lemma 2.4.
\endroster
\QED

All of the steps here can be formalized in primitive recursive arithmetic,
so Theorem~4.4 is a theorem of PRA.

\bigpagebreak
\head 5. Embedding algebras \endhead

In this section, we consider algebras of increasing functions from $\natnum$ to
$\natnum$ which imitate the behavior of the algebra of elementary embeddings
from Laver~\cite{\LavLDLFA} when restricted to the set of critical points.  The
existence of such algebras will turn out to be equivalent to the properties
in Theorem~4.4. Moreover, this equivalence can be proved (and formulated)
in primitive recursive arithmetic.

Let $\id$ be the identity function on~$\natnum$.  If
$f\colon\natnum\to\natnum$ is strictly
increasing and different from~$\id$, let $\crit f$ be the least $n$ such
that $f(n)>n$ (the {\it critical point\/} of $f$).

\definition{Definition 5.1} An {\it embedding algebra}
is a structure~$(A,{\cdot})$ where $A$~is a collection of
strictly increasing functions from~$\natnum$ to~$\natnum$, ${\cdot}$~is a
left-distributive binary operation on $A$, 
and, for every $a, b \in A$ with $b\neq\id$, $\crit{a\cdot b}=a(\crit b)$.
\enddefinition

As usual, we will often write~$ab$ instead of~$a\cdot b$.  The set $A$
need not contain the identity function, but, if it does not, one can
extend the operation $\cdot$ to $A \cup \{\id\}$ in the obvious way:
$a\cdot\id = \id$, $\id\cdot a = a$.

An embedding algebra~$A$ is {\it nontrivial\/} if it has an element
other than $\id$.  Note that the set of non-identity elements of~$A$ is
closed under~$\cdot$: if $b$ has a critical point, so does $a \cdot b$.
Also, if $A$~is nontrivial, then $A$~has infinitely many critical
points: if $n = \crit a$, then $a(n) = \crit{aa}$ and $a(n) > n$.

The main goal of the next three sections will be to prove the following
theorem.

\proclaim{Theorem 5.2} The statement ``There exists a
nontrivial embedding algebra''
is equivalent to the statement ``$A_\infty$ is free''. \endproclaim

When proving that ``$A_\infty$ is free'' implies the existence of an
embedding algebra, we shall see that there is a natural way of associating
increasing functions from $\natnum$ to $\natnum$ with words in $W_\cala$.
However, it is not easy to prove that inequivalent words yield distinct
functions. Here we shall rely on Theorem~2.14, but first we
have to develop techniques to `miniaturize' Laver's proof. This will be
done in Section 6. In order to develop the necessary machinery, we first
define a different kind of `embedding algebra.'
The new definition will include much of Laver's machinery
explicitly; the resulting
structure will be much less concrete but more amenable to algebraic
manipulation.

\definition{Definition 5.3}
A {\it two-sorted embedding algebra} consists of a nonempty
set $\E$ (the `embeddings,' for which we will use variables $a,b,\dotsc$) and
a nonempty set $\O$ (the `ordinals,' for which we will use variables
$\alpha,\beta,\dotsc$), together with binary operations $\cdot$ and $\circ$ on
$\E$, a binary relation $\le$ on $\O$, a constant $\id \in \E$, an application
operation $a,\beta \mapsto a(\beta)$ (which will often be written
without parentheses) from $\E \times \O$ to $\O$, a function
$\critfunc\colon \E{-}\{\id\} \to \O$, and a ternary relation ${\equiv}
\subseteq \E \times \O \times \E$, satisfying the following axioms:
\roster
\item"$\bullet$" The relation $\le$ is a linear ordering of $\O$.
\nopagebreak
\item"$\bullet$" Embeddings are strictly increasing monotone functions:
$$\beta < \gamma \quad \text{implies} \quad a\beta < a\gamma,
\qquad \text{and} \qquad a\beta \ge \beta.$$
\item"$\bullet$" For all $a \ne \id$, $a(\crit a) > \crit a$.
\item"$\bullet$" The operation $\circ$ represents composition:
$(a\circ b)\gamma = a(b\gamma)$.
\item"$\bullet$" The constant $\id$ represents the identity:
$$\id(\gamma) = \gamma, \qquad a \cdot \id = \id, \qquad \text{and} \qquad
\id \cdot a = a \circ \id = \id \circ a = a.$$
\item"$\bullet$" The axioms (\Lavlaws) hold.
\item"$\bullet$" For each $\gamma$, $\tLequiv\gamma/$ is an equivalence
relation
on $\E$ which respects $\cdot$ and $\circ$
(i.e., if $a \tLequiv\gamma/ a'$ and $b \tLequiv\gamma/ b'$,
then  $a\cdot b \tLequiv\gamma/ a'\cdot b'$
and $a\circ b \tLequiv\gamma/ a'\circ b'$).
\item"$\bullet$" If $\gamma \le \delta$ and
$a \tLequiv\delta/ b$, then $a \tLequiv\gamma/ b$.
\item"$\bullet$" If $a \tLequiv\gamma/ b$ and $a\delta<\gamma$, then
$a\delta = b\delta$.
\item"$\bullet$" For any $a \ne \id$, $a \tLequiv\crit a/ \id$.
\item"$\bullet$" Coherence: $a \tLequiv\gamma/ b$ implies
$ca \tLequiv c\gamma/ cb$.
\endroster
\enddefinition

It follows from these axioms that the operation $\cdot$ distributes over
itself and application:
$$a(bc)=ab(ac) \qquad\text{and}\qquad a(b\gamma)=ab(a\gamma).$$
A few more properties also follow easily:

\proclaim{Proposition 5.4} In a two-sorted embedding algebra, if $a$
and $b$ are embeddings different from $\id$, then:
\roster
\item $\crit a$ is the least ordinal moved by $a$;
\item $\crit{ab} = a(\crit b)$; and
\item $\crit{a\circ b} = \min(\crit a,
\crit b)$.
\endroster
\endproclaim

\demo{Proof} It is given that $a(\crit a) > \crit a$; if $\beta < \crit a$,
then the fact that $\id \tLequiv\crit a/ a$ implies that $\beta = \id(\beta)
= a(\beta)$.  Since $\id \tLequiv\crit b/ b$, coherence gives
$\id = a\cdot\id \tLequiv a(\crit b)/ ab$, so $ab$ does not move
any ordinal less than $a(\crit b)$; but it moves $a(\crit b)$ to
$ab(a(\crit b)) = a(b(\crit b)) > a(\crit b)$, so we must have
$\crit{ab} = a(\crit b)$.  For \therosteritem3, let $\gamma = \min(\crit a,
\crit b)$.  Then, since $\tLequiv\gamma/$ respects $\circ$, we have
$\id \tLequiv\gamma/ a \circ b$, while $(a\circ b)\gamma =
a(b\gamma) \ge \max(a\gamma,b\gamma) > \gamma$, so $\gamma$
is the least ordinal moved by $a \circ b$. \QED

It is easy to verify that all of the axioms in Definition~5.3
are preserved when one moves
to a substructure (replacing $\E$ and $\O$ with smaller sets closed
under the operations, and restricting the operations and relations
accordingly).  In particular, if one keeps the same $\E$ but replaces~$\O$
with the range of the function $\critfunc$ (assuming that $\E \ne
\{\id\}$), then Proposition~5.4(2) implies that
the new sets are closed under the operations, 
so one obtains a new
two-sorted embedding algebra in which every ordinal is a critical point.

If desired, one can restrict $\E$ to the embeddings obtained from a
single embedding $j\ne\id$ using $\cdot$ and $\circ$, along with~$\id$;
this gives a two-sorted embedding algebra generated by a single embedding.
From now on, we will call a two-sorted embedding algebra {\it monogenic}
if its non-identity embeddings are generated from a single non-identity
embedding via $\cdot$ and~$\circ$.  Similarly, an embedding algebra is
{\it monogenic} if it is generated from a single non-identity embedding
via~$\cdot$; any nontrivial embedding algebra has monogenic
subalgebras.  Note that a monogenic embedding algebra does not contain
the identity function.

The results of Laver~\cite{\LavLDLFA} show that one can make the set of all
elementary embeddings from $V_{\lambda}$ to itself into a two-sorted
embedding algebra by letting $\O$ be the set of limit ordinals less
than $\lambda$ and defining $\tLequiv\gamma/$
to be $\overset{\gamma}\to=$ (as defined in
Laver~\cite{\LavLDLFA}, Section 2).  
We now want to show that just the
simple properties of embedding algebras suffice to construct the more
elaborate apparatus of a two-sorted embedding algebra.

\proclaim{Proposition~5.5} If a nontrivial embedding algebra exists,
then there exists a two-sorted embedding algebra in which
the ordinals have order type~$\omega$. \endproclaim

\demo{Proof}
Let such an embedding algebra be given; we will construct a
two-sorted embedding algebra.  The ordinal set~$\O$ will be the set
of critical points from the given algebra; this is an infinite subset
of~$\natnum$, so it has order type~$\omega$.  The embedding set
and the operations and relations will be built up in several steps.

To start with, let $\E_1$ be the set of non-identity embeddings
in the given algebra.  As noted before, this set is closed under~$\cdot$.
Now the following properties are true of $\E_1$ and~$\O$:
\roster
\item"$\bullet$" The left distributive law holds.
\item"$\bullet$" $\beta < \gamma$ implies $a\beta < a\gamma$.
\item"$\bullet$" $a(\gamma) \ge \gamma$.
\item"$\bullet$" $a(\crit a) > \crit a$.
\item"$\bullet$" $a(\gamma) = \gamma$ for $\gamma < \crit a$.
\item"$\bullet$" $\crit{ab} = a(\crit b)$.
\endroster
We also have the property
\roster
\item"$\bullet$" $ab(a\gamma) = a(b\gamma)$,
\endroster
since every ordinal~$\gamma$ is a critical point and
$$ab(a(\crit c)) = ab(\crit{ac})= \crit{ab(ac)} = \crit{a(bc)}
= a(\crit{bc}) = a(b(\crit c)).$$

Now use the construction from Proposition~2.10
to extend and expand $(\E_1,{\cdot})$ to an algebra
$(\E_2,{\cdot},{\circ})$ satisfying Laver's laws~(\Lavlaws).
The application operation on these new embeddings
is defined naturally: each embedding~$a$ is a formal composition
$(a_1 \circ \dots \circ a_n)$ of members of~$\E_1$, and we let
$a(\gamma) = a_1(a_2(\dots a_n(\gamma)\dots))$.
We have $a_ia_{i+1}(a_i(\delta)) = a_i(a_{i+1}(\delta))$
for any $\delta$, so replacing $a_i \circ a_{i+1}$
with $a_ia_{i+1} \circ a_i$ in the formal composition
does not change the resulting value of $a(\gamma)$; since
formal compositions were identified only when one could transform
one into the other by such replacements and/or the reverse,
the value $a(\gamma)$ is well-defined.
Also, let $\crit a$ be the minimum of $\crit{a_1},\dots,\crit{a_n}$;
this is the least~$\gamma$ such that $a(\gamma) > \gamma$,
so it also does not depend on the expression for~$a$.
Then we have:
\roster
\item"$\bullet$" (\Lavlaws) holds.
\item"$\bullet$" $(a\circ b)\gamma = a(b\gamma)$.
\item"$\bullet$" $\crit{a\circ b} = \min(\crit a, \crit b)$.
\endroster
And the properties listed before hold for~$\E_2$ as well.

Let~$\E$ be $\E_2 \cup \{\id\}$, where $\id$ is a new embedding
for which $\crit{\id}$ is not defined but the other operations
are defined by:
\roster
\item"$\bullet$" $\id(\gamma) = \gamma$, $a \cdot \id = \id$, and
$\id \cdot a = a \circ \id = \id \circ a = a$.
\endroster
Again the previous properties continue to hold.  Now it only
remains to define $a \tLequiv\gamma/ b$ so that the
rest of the axioms in Definition~5.3 hold.

\proclaim{Lemma 5.6} Assume the facts listed above.
Let $a,b_1,\dots,b_k$ be embeddings, where $k \ge 0$,
and let $\gamma$ be an ordinal.

{\rm (i)} If $\crit a > b_1b_2\dotsm b_k\gamma$, then
$a b_1b_2\dotsm b_k\gamma =  b_1b_2\dotsm b_k\gamma$.

{\rm (ii)} If $\crit a > ab_1b_2\dotsm b_k\gamma$, then
$a b_1b_2\dotsm b_k\gamma =  b_1b_2\dotsm b_k\gamma$.
\endproclaim

\demo{Proof} These are both proved by induction on $k$
(simultaneously for all embeddings).
Let us write (i${}_m$) for the case $k = m$ of~(i),
and similarly for (ii).  Note that the hypotheses of (i)
and (ii) each imply that $\crit a > \gamma$.

(i${}_0$): This just says that $a$ does not move any ordinal below its
critical point.

(i${}_1$): $ab_1\gamma = ab_1(a\gamma) = a(b_1\gamma) = b_1\gamma$.

(i${}_k$) for $k \ge 2$: Let $s = ab_1a$.  Note that $s(ab_1b_2) =
ab_1a(ab_1b_2) = ab_1(ab_2) = a(b_1b_2)$.  Also note that $\crit s =
ab_1(\crit a) \ge \crit a$; similarly, $\crit{ws} \ge \crit a$ for
any $w$.  In particular,
$$\crit{b_1b_2 \dotsm b_{k-1}s} \ge \crit a > b_1b_2 \dotsm b_k \gamma,$$
so (i${}_1$) gives
$$b_1b_2\dotsm b_{k-1}(sb_k)\gamma = b_1b_2\dotsm b_{k-1}s(b_1b_2\dotsm
b_{k-1}b_k)\gamma = b_1b_2\dotsm b_{k-1}b_k\gamma.$$
We now have
$$\crit{b_1b_2 \dotsm b_{k-2}s} \ge \crit a > b_1b_2 \dotsm b_{k-1}(sb_k)
\gamma,$$
so, if $k > 2$, we can apply (i${}_2$) to get
$$\multline
b_1b_2\dotsm b_{k-2}(sb_{k-1})(sb_k)\gamma = b_1b_2\dotsm
b_{k-2}s(b_1b_2\dotsm b_{k-2}b_{k-1})(sb_k)\gamma \\
= b_1b_2\dotsm b_{k-1}(sb_k)\gamma =
b_1b_2\dotsm b_{k-1}b_k\gamma.
\endmultline$$
We can now apply (i${}_3$) to $b_1b_2\dotsm b_{k-3}s$, and so on
all the way to (i${}_{k-2}$), to get
$$b_1b_2(sb_3)(sb_4)\dotsm(sb_k)\gamma = b_1b_2\dotsm b_k\gamma.$$
Now we have
$$\alignat2
s(ab_1b_2\dotsm b_k\gamma)
&= s(ab_1b_2)(sb_3)(sb_4)\dotsm (sb_k)(s\gamma) &\qquad\qquad&\\
&= a(b_1b_2)(sb_3)(sb_4)\dotsm (sb_k)\gamma &&\\
&=b_1b_2(sb_3)(sb_4)\dotsm (sb_k)\gamma &&\text{by (i${}_{k-1}$)}\\
&= b_1b_2\dotsm b_k\gamma &&\\
&= s(b_1b_2\dotsm b_k\gamma). &&
\endalignat$$
Since $s$ maps distinct ordinals to distinct ordinals, we get
$a b_1b_2\dotsm b_k\gamma =  b_1b_2\dotsm b_k\gamma$.

(ii${}_0$): We have $\crit a > a\gamma \ge \gamma$, so (i${}_0$) applies.

(ii${}_1$): $ab_1\gamma = ab_1(a\gamma) = a(b_1\gamma) \ge b_1\gamma$,
so $\crit a > b_1\gamma$, so (i${}_1$) applies.

(ii${}_k$) for $k \ge 2$: Again let $s = ab_1a$.  We now have
$$\crit{ws} \ge \crit s \ge \crit a > a b_1b_2\dotsm b_k\gamma$$
for any $w$.  This gives
$$\alignat2
ab_1b_2\dotsm b_k\gamma
&= s(ab_1b_2\dotsm b_k\gamma) &\qquad\qquad&\\
&= s(ab_1b_2)(sb_3)(sb_4)\dotsm (sb_k)(s\gamma) &&\\
&= a(b_1b_2)(sb_3)(sb_4)\dotsm (sb_k)\gamma &&\\
&= b_1b_2(sb_3)(sb_4)\dotsm (sb_k)\gamma &&\text{by (ii${}_{k-1}$)}\\
&= b_1b_2s(b_1b_2b_3)(sb_4)\dotsm (sb_k)\gamma &&\\
&= b_1b_2b_3(sb_4)\dotsm (sb_k)\gamma &&\text{by (ii${}_{k-2}$)}\\
&= b_1b_2b_3s(b_1b_2b_3b_4)(sb_5)\dotsm (sb_k)\gamma &&\\
&= b_1b_2b_3b_4(sb_5)\dotsm (sb_k)\gamma &&\text{by (ii${}_{k-3}$)}\\
&= \dotsb &&\\
&= b_1b_2\dotsm b_k\gamma, &&\text{by (ii${}_1$)}
\endalignat$$
as desired.
\QED

Define the preliminary relation $\Lequiv\gamma/$ between embeddings
as follows:
$a \Lequiv \gamma/ b$ if, for each $k \ge 0$ and all embeddings
$c_1,\dots,c_k$,
$$ac_1\dotsm c_k \Lrestrict \gamma = bc_1\dotsm c_k \Lrestrict \gamma,$$
where $a \Lrestrict \gamma$ is
$a \restrict \{\beta\colon a(\beta) < \gamma\}$.
In other words, $a \Lequiv \gamma/ b$ iff, for any $\delta$, if either
$ac_1\dotsm c_k\delta$ or $bc_1\dotsm c_k\delta$
is less than $\gamma$, then $ac_1\dotsm c_k\delta = bc_1\dotsm c_k\delta$.
This is easily seen to be an equivalence relation, and
Lemma~5.6 just states that $a \Lequiv\crit a/ \id$.

We can now define the final desired relation $\tLequiv\gamma/$ by:
$a \tLequiv\gamma/ b$ iff
$ra \Lequiv r\gamma/ rb$ for all embeddings $r$ (including
$r = \id$).  This is also an equivalence relation.  Since
$\crit{ra} = r(\crit a)$, we have $a \tLequiv\crit a/\id$.

If $a \tLequiv\gamma/ b$, then $(r\circ c)a \Lequiv(r\circ c)\gamma/
(r\circ c)b$ for any $r$, so $r(ca) \Lequiv r(c\gamma)/ r(cb)$;
hence, $ca \tLequiv c\gamma/ cb$.

Easily, if $\gamma \le \delta$, then $a \Lequiv\delta/ b$ implies
$a \Lequiv\gamma/ b$, and the same holds for $\tLequiv/$.

It follows immediately from the definitions of $\tLequiv\gamma/$
(with $r = \id$) and $\Lequiv\gamma/$
(with $k=0$) that, if $a \tLequiv\gamma/ b$ and $a\delta < \gamma$,
then $a\delta = b\delta$.

If $a \tLequiv\gamma/ a'$ and $b \tLequiv\gamma/ b'$, then
we have already shown that $ab \tLequiv a\gamma/ ab'$,
so $ab \tLequiv \gamma/ ab'$.  Also,
$r(ab)c_1\dotsm c_k = ra(rb)c_1\dotsm c_k$ and
$r(a'b)c_1\dotsm c_k = ra'(rb)c_1\dotsm c_k$, so from
$a \tLequiv\gamma/ a'$ we get $ab \tLequiv \gamma/ a'b$.
Similarly, we get $(a\circ b) \tLequiv\gamma/ (a' \circ b)$
since $r(a\circ b) c_1\dotsm c_k = ra(rbc_1)c_2\dotsm c_k$
and the same for $a'$.  (For the case $k = 0$, note that,
if $r(a \circ b)\delta < r\gamma$, then $ra(rb\delta) < r\gamma$,
so $ra(rb\delta) = ra'(rb\delta)$, so $r(a\circ b)\delta =
r(a' \circ b)\delta$.)  Now, using the formulas
$a \circ b = ab \circ a$ and $a \circ b' = ab' \circ a$, we get
$(a \circ b) \tLequiv\gamma/ (a \circ b')$.
So the equivalence relation $\tLequiv\gamma/$ respects application
and composition of embeddings.

Therefore, we have a two-sorted embedding algebra.
\QED

If the original embedding algebra satisfies the property
$ab(a(n))=a(b(n))$ for all embeddings $a,b$ and natural numbers~$n$,
then one can let~$\O$ be the entire set~$\natnum$, rather
than just the critical points, and the construction
will work as before.  As a result, one sees
that the two-sorted embedding algebra includes an `isomorphic'
copy of the original embedding algebra, expressed in two-sorted form.
[In order to see that deleting~$\id$ and reinserting it later
does not cause a problem, we must show that the new formulas for
multiplying by~$\id$ match the old ones.  In other words, we
must see that, if the original embedding algebra contained~$\id$,
then it satisfied $\id \cdot a = a$ and $a \cdot \id = \id$.
To see this, use the property above to get, for all~$n$,
$$(\id \cdot a)(n) = (\id \cdot a)(\id(n)) = \id(a(n)) = a(n)$$ and
$$(a \cdot \id)(a(n)) = a(\id(n)) = a(n) = \id(a(n)).$$
So $\id \cdot a = a$, and $a \cdot \id$ agrees with~$\id$ at all
numbers of the form~$a(n)$; but the only strictly increasing
function from~$\natnum$ to~$\natnum$ which agrees with~$\id$
at infinitely many places is~$\id$.]

It is easy to see that, if the original embedding algebra is monogenic,
then so is the two-sorted embedding algebra constructed above.

We conclude this section with a proposition about two-sorted embedding
algebras which
is a substitute for Kunen's theorem about elementary
embeddings.

For any non-identity embedding $a$, the sequence $\crit a, a(\crit a),
a(a(\crit a)), \dotsc$ is a strictly increasing sequence of ordinals,
called the {\it critical sequence} of $a$.

\proclaim{Proposition 5.7} In any monogenic two-sorted embedding algebra,
if $a\ne\id$ is an embedding, then the critical sequence of $a$ is
cofinal in the set of critical points (the range of~$\critfunc$).
Also, $a(\gamma) > \gamma$ for any critical point $\gamma \ge \crit a$.
\endproclaim

\demo{Proof} All members of the critical sequence are critical points
(of the embeddings $a$,~$aa$, $a(aa)$, $a(a(aa))$, etc.).
Let $j$ be a non-identity embedding which generates the algebra, and let
$\langle \kappa_n \colon n \in \natnum\rangle$ be the critical sequence of $j$.
We recall Lemma 2.15. It was stated for elementary embeddings, but the proof
clearly works in the present context as well. Thus every
$a$ must move some ordinal $\kappa_n$, and hence
$\crit a \le \kappa_n$; this shows that the critical sequence of $j$
is cofinal in the critical points.  To complete the proof
of the first claim, we now
show by induction on expressions in $j$ that, if $a \ne \id$
and $\langle \alpha_n \colon n \in \natnum\rangle$ is the critical sequence
of $a$, then $\alpha_n \ge \kappa_n$ for all $n$.  This is again
trivial for $a=j$.  Suppose it is true for $b$ and $c$, with critical
sequences $\langle \beta_n\colon n \in \natnum\rangle$ and
$\langle \gamma_n\colon n \in \natnum\rangle$ respectively.  If
$a = bc$, then induction gives $\alpha_n = b\gamma_n$ for all $n$,
so $\alpha_n = b\gamma_n \ge \gamma_n \ge \kappa_n$.  If $a =
b \circ c$, then $\alpha_0$ is either $\beta_0$ or $\gamma_0$.
In the former case, the fact that $\alpha_{n+1} = b(c\alpha_n) \ge
b\alpha_n$ gives $\alpha_n \ge \beta_n$ for all $n$; similarly,
in the latter case, we have $\alpha_n \ge \gamma_n$ for all $n$.
In either case, we get $\alpha_n \ge \kappa_n$, as desired.

Now, if $\gamma \ge \crit a$ is a critical point, then $\gamma \ge \alpha_0$
and $\gamma < \alpha_m$ for some~$m$, so there is an~$n$ such that
$\alpha_n \le \gamma < \alpha_{n+1}$.  This gives
$a\gamma \ge a\alpha_n = \alpha_{n+1} > \gamma$. \QED

%

\bigpagebreak
\head 6. Extended two-sorted embedding algebras \endhead

In order to prove Theorem~5.2, we will need to perform a number of the
arguments of Laver~\cite{\LavOAEER} in the context of two-sorted
embedding algebras.
This is straightforward for arguments involving only the operations
which are built into these algebras, but some arguments use additional
features of elementary embeddings.  In particular, a few arguments
use ordinals of the form $\appless a\gamma$, defined to be the least ordinal
greater than $a(\beta)$ for all $\beta < \gamma$.  In this section,
we will define an extended algebra which includes this operation and
show that such algebras can be constructed from ordinary two-sorted
embedding algebras; this will allow us to use this new operation
to prove facts about the original algebra.

\definition{Definition 6.1} An {\it extended two-sorted embedding algebra}
is a two-sorted embedding algebra (with embedding
set $\E$ and ordinal set $\O$), together with two new operations,
a cofinality function ${\cf}\colon\O\to\O$ and a mapping from
$\E\times\O$ to $\O$ for which we use the notation $a,\gamma \mapsto
\appless a\gamma$, satisfying the following additional axioms:
$$\allowdisplaybreaks\gather
a(\appless b\gamma) = \appless{ab}{a\gamma}; \\
\appless a{\appless b\gamma} = \appless{(a \circ b)}\gamma; \\
\appless a\gamma \le a\gamma; \\
\text{if}\quad \gamma<\delta, \quad\text{then}\quad a\gamma<\appless a\delta;\\
\text{if}\quad a \tLequiv\gamma/ b \quad\text{and}\quad \appless a\delta \le
   \gamma, \quad\text{then}\quad \appless a\delta = \appless b\delta; \\
\cf(\crit a) = \crit a; \\
\cf(\appless a\gamma) = \cf\gamma; \\
\cf(a\gamma) = a(\cf\gamma); \\
\cf\gamma \le \gamma; \\
\text{if}\quad a(\cf\gamma) = \cf\gamma, \quad\text{then}\quad \appless a\gamma
   = a\gamma.
\endgather $$
\enddefinition

The last two of these axioms are not used in this paper, but they might be
useful for later applications.  On the other hand, there are a few
facts that are used in this paper but not given above, because they
can be deduced from the axioms.

\proclaim{Proposition 6.2} In an extended two-sorted embedding algebra:
\roster
\item $\appless a\gamma \ge \gamma$;
\item $\appless \id\gamma = \gamma$;
\item $\appless a\gamma = \gamma$ for $\gamma \le \crit a$; and
\item if $a(\cf\gamma) > \cf\gamma$, then $\appless a\gamma < a\gamma$.
\endroster \endproclaim

\demo{Proof} For all $\delta < \gamma$, we have $\appless a\gamma > a\delta
\ge \delta$; hence, \therosteritem1 holds.  This and $\appless \id\gamma \le
\id(\gamma)$ give~\therosteritem2; we then get~\therosteritem3 because
$a \tLequiv\gamma/ \id$. For~\therosteritem4, we have $\appless a\gamma \le
a\gamma$, and equality cannot hold because $\appless a\gamma$ and $a\gamma$
have different cofinalities. \QED

Again it is not hard to verify that the axioms for an extended two-sorted
embedding algebra hold in the case where $\E$ is a set of elementary
embeddings on $V_\lambda$ and $\O$ is the collection of limit ordinals less
than $\lambda$~\cite{\LavOAEER}.  Also, any subalgebra of an extended
two-sorted
embedding algebra is also an extended two-sorted embedding algebra; in
particular, if we keep the same set of embeddings but restrict the ordinals to
those of the form $\appless a{\crit b}$, we get an algebra in which all
ordinals
have this form.  (Proposition~6.2\therosteritem3 gives
$\crit a = \appless a{\crit
a}$, so all critical points are in this set of ordinals; now the axioms easily
imply that this set of ordinals is closed under all of the algebra
operations.)

We now state the main result of this section.

\proclaim{Theorem 6.3} Suppose that we are given a two-sorted embedding
algebra, in which every ordinal is a critical point.
Then the algebra can be extended to a new two-sorted embedding algebra
with the same embedding set, on which the required additional operations
can be defined so as to give an extended two-sorted embedding algebra.
\endproclaim

The proof of this theorem will use the following two lemmas about
two-sorted embedding algebras.

\proclaim{Lemma 6.4} In any two-sorted embedding algebra, if\/ $\gamma =
\crit c$, then:
\roster
\item"(a)" $cc(ca\gamma) < c(ca\gamma)$;
\item"(b)" $ca\gamma$ is not in the range of $c$.
\endroster
\endproclaim

\demo{Proof} For (a), note that $\crit{cc} = c\gamma > \gamma$, so
$cc\gamma = \gamma$; hence,
$$c(ca\gamma) = cc(ca)(c\gamma) > cc(ca)\gamma = cc(ca)(cc\gamma)
= cc(ca\gamma).$$
On the other hand, an element $\delta$ of the range of $c$ cannot satisfy
$cc\delta<c\delta$; if $\delta = c\beta$, then
$c\delta = cc(c\beta) = cc\delta$.  Therefore, (b) holds. \QED

\proclaim{Lemma 6.5} In any two-sorted embedding algebra, if\/
$\crit r = \crit s = \kappa$, then $r\lambda < ra\kappa$ implies
$s\lambda < sa\kappa$. \endproclaim

\demo{Proof} Assume $r\lambda < ra\kappa$.  Note that $\crit{rs} =
r\kappa > \kappa$, so $rs\kappa = \kappa$; this gives
$$r(s\lambda) = rs(r\lambda) < rs(ra\kappa) = rs(ra)(rs\kappa)
= rs(ra)\kappa < rs(ra)(r\kappa) = r(sa\kappa).$$
Since $r$ gives an increasing function on the ordinals, we must have
$s\lambda < sa\kappa$. \QED

\demo{Proof of Theorem 6.3}
Fix a two-sorted embedding algebra.  Let $\E$ and
$\O$ be its embedding set and ordinal set, respectively, and assume
that the range of $\critfunc$ is all of $\O$.  We must extend $\O$ to a larger
collection of ordinals on which the operation $\appless a\gamma$ can
be suitably defined.  The remarks following Proposition~6.2 indicate
that this new set of ordinals need only contain the ordinals
$\appless a{\crit b}$ for $a,b \in \E$.  The main step will be to define
the linear ordering properly for such ordinals; it turns out that
the properties of an extended two-sorted embedding algebra determine
this ordering completely.

\proclaim{Lemma 6.6} In an extended two-sorted embedding algebra, if\/
$\gamma = \crit c$ and $\delta$ is any ordinal, then
$$\appless a\gamma \le \delta \iff ca\gamma < c\delta.$$
\endproclaim

\demo{Proof} If $\appless a\gamma \le \delta$, then the fact that $c\gamma >
\gamma$ gives $$ca\gamma < \appless{ca}{c\gamma}= c(\appless a\gamma) \le
c\delta.$$  On the other hand, if $\delta < \appless a\gamma$, then we can
use $\appless c\gamma=\gamma$ to get
$$c\delta < \appless c{\appless a\gamma} = \appless{(c\circ a)}\gamma =
\appless{(ca \circ c)}\gamma = \appless{ca}{\appless c\gamma}=
\appless{ca}\gamma
\le ca\gamma.$$ \QED

It follows that, if $\gamma = \crit c$ and $\delta = \crit d$, then
$$\align \appless a\gamma \le \appless b\delta
&\iff ca\gamma < c(\appless b\delta) \\
&\iff \appless {cb}{c\delta} \not\le ca\gamma \\
&\iff cd(cb)(c\delta) \not< cd(ca\gamma) \\
&\iff cd(ca\gamma) \le c(db\delta).
\endalign$$
This tells us how to start the construction from the given two-sorted
embedding algebra.

We want to define a binary relation $R$ on $\E\times\O$ as follows:
$$(a,\gamma)R(b,\delta) \iff cd(ca\gamma) \le c(db\delta),$$ where $c$ and
$d$ are chosen so that $\crit c = \gamma$ and $\crit d = \delta$.  Such $c$
and $d$ do exist because every element of $\O$ is a critical point; we must
now see that the definition of $R$ does not depend on which $c$ and $d$ are
chosen.  If $c'$ also has critical point $\gamma$, then Lemma~6.5 gives
$$c(d\circ a)\gamma \le c(db\delta) \iff c'(d\circ a)\gamma \le
c'(db\delta),$$ so $cd(ca\gamma) \le c(db\delta)$ iff $c'd(c'a\gamma)
\le c'(db\delta)$.  Also, if $d'$ is another embedding with critical point
$\delta$, then $\crit{cd} = \crit{cd'} = c\delta$, so Lemma~6.5 gives
$$cd(ca\gamma) < cd(cb)(c\delta) \iff cd'(ca\gamma) <
cd'(cb)(c\delta).$$
Note that, by Lemma~6.4(b), $c(d\circ a)\gamma \le c(db\delta)$ is
equivalent to $c(d\circ a)\gamma < c(db\delta)$, so
$(a,\gamma)R(b,\delta)$ iff $cd(ca\gamma) < c(db\delta)$.
Therefore, $R$ is well-defined.

Lemma~6.4(a) implies that $R$ is reflexive.  We will now show that $R$ is
transitive.  Suppose $(a,\rho)R(b,\sigma)$ and $(b,\sigma)R(c,\tau)$;
fix embeddings $r,s,t$ with critical points $\rho,\sigma,\tau$,
respectively.  We then have $rs(ra\rho) \le r(sb\sigma)$ and
$st(sb\sigma) \le s(tc\tau)$, so
$$\align
rs(rt(ra\rho)) &= rs(rt)(rs(ra\rho)) \\
&\le rs(rt)(r(sb\sigma)) \\
&= r(st(sb\sigma)) \\
&\le r(s(tc\tau)) \\
&= rs(r(tc\tau)),
\endalign$$
so $rt(ra\rho) \le r(tc\tau)$, so $(a,\rho)R(c,\tau)$.
The same proof using $>$ instead of $\le$ shows that the negation of $R$
is also transitive.

We now know that $R$ is a preorder; if we define the relation $\sim$ on
$\E\times\O$ by
$$(a,\gamma)\sim(b,\delta) \iff (a,\gamma)R(b,\delta) \text{ and }
(b,\delta)R(a,\gamma),$$
then $\sim$ is an equivalence relation on $\E\times\O$ and $R$ induces
a partial order on the set of equivalence classes.  Let $\O^*$ be the
set of equivalence classes; we will write $[a,\gamma]$ for the
equivalence class of $(a,\gamma)$.  Let $\le^*$ be the partial ordering
induced by $R$ on $\O^*$.  We then have
$$\align [a,\gamma] \le^* [b,\delta]
&\iff cd(ca\gamma) \le c(db\delta) \\
&\iff cd(ca\gamma) < c(db\delta),
\endalign$$
where $\crit c = \gamma$ and $\crit d = \delta$.  The fact that the negation
of $R$ is transitive implies that any two elements of $\E\times\O$ are
$R$-comparable (if $xRy$ and $yRx$ were both false, then $xRx$ would be false,
contradicting reflexivity), so $\le^*$ is a linear ordering of $\O^*$.

The various distributive laws imply that, for any $e \in \E$, we have
$(a,\gamma) R (b,\delta)$ if and only if $(ea,e\gamma) R (eb,e\delta)$.
Therefore, $e$ induces a mapping from $\O^*$ to $\O^*$ via the formula
$e[a,\gamma] = [ea,e\gamma]$, and this mapping is strictly increasing.
Also, we clearly have $(e\circ e')[a,\gamma] = e(e'[a,\gamma])$.

The element $[a,\gamma]$ of $\O^*$ is meant to represent $\appless a\gamma$
in an
extended algebra.  For this to extend the original algebra, we need an element
$H(\gamma)$ of $\O^*$ to correspond to each $\gamma \in \O$.  This element
will turn out to be $[c,\gamma]$, where $c$ is any embedding with critical
point $\gamma$.  In order to see that this is well-defined and gives the
proper ordering on the representatives in $\O^*$, we need the following
result.

\proclaim{Lemma 6.7} If\/ $\crit c = \gamma$ and\/ $\crit d = \delta$, then\/
$[c,\gamma] \le^* [d,\delta]$ if and only if\/ $\gamma \le \delta$.
\endproclaim

\demo{Proof} By definition, $[c,\gamma] \le^* [d,\delta]$ if and only if
$cd(cc\gamma) \le c(dd\delta)$.  But $\crit{cc} > \gamma$ and
$\crit{dd} > \delta$, so this is equivalent to $cd\gamma \le c\delta$.
Now, if $\gamma = \delta$, then $\crit{cd} = c\gamma > \gamma$, so
$cd\gamma = \gamma = \delta \le c\delta$.  If $\gamma > \delta$,
then $cd\gamma \ge \gamma > \delta = c\delta$, so
$[c,\gamma] \not\le^* [d,\delta]$, so $[c,\gamma] >^* [d,\delta]$.
Symmetrically, if $\gamma < \delta$, then $[c,\gamma] <^* [d,\delta]$. \QED

So the correspondence between $\gamma$ and $[c,\gamma]$ gives
an order-preserving map ${H\colon\O\to\O^*}$.  This lets us define
the new critical point map $\critfunc^* \colon \E \to \O^*$ by
the formula $\critfunc^*(c) = H(\crit c) = [c,\crit c]$.

We next verify that the embedding maps $\gamma^* \mapsto e\gamma^*$ satisfy
$e[a,\gamma] \ge^* [a,\gamma]$.  We must show that $c(ec)(ca\gamma) \le
c(ec(ea)(e\gamma))$, where $\crit c = \gamma$; to see this, note that
$$\align
c(ec)(ca\gamma) &= ce(cc)(ca\gamma) \\
&\le ce(cc)(ce(ca\gamma)) \\
&= ce(cc(ca\gamma)) \\
&< ce(c(ca\gamma)) \qquad\qquad\text{by Lemma 6.4(a)} \\
&= c(e(ca\gamma)) \\
&= c(ec(ea)(e\gamma)).
\endalign$$

Clearly $e(\critfunc^*(a)) = \critfunc^*(ea)$; since $\crit{aa} =
a(\crit a) > \crit a$, this gives $a(\critfunc^*(a)) = \critfunc^*(aa) >
\critfunc^*(a)$.

Next, we define the new ternary relation $\tLequiv*/$ as follows:
$a \tLequiv*\,\gamma^*/ b$ iff $a \tLequiv\delta/ b$ for some
$\delta\in\O$ such that $\gamma^* \le^* H(\delta)$.  In other words, $a$ agrees
with $b$ up to some new ordinal iff $a$ agrees with $b$ up to some
old ordinal at least as high.
Using this definition, it is easy to deduce all of the axioms about
$\tLequiv*/$ from the corresponding axioms about $\tLequiv/$, except
for the axiom ``if $a \tLequiv*\,\gamma^*/ b$ and $a\delta^* < \gamma^*$,
then $a\delta^* = b\delta^*$''; this one will require more work.

If $\crit d = \delta$, then $H(\delta) = [d,\delta] \le^* [a,\delta]$ for any
$a$, because $dd(dd\delta) = \delta \le d(da\delta)$.

\proclaim{Lemma 6.8} If\/ $\crit c = \gamma$, then\/ $[a,\gamma] \le^*
H(\delta)$ if and only if $ca\gamma < c\delta$.  \endproclaim

\demo{Proof} Fix $d$ with critical point $\delta$; then $[a,\gamma]
\le^* [d,\delta]$ is equivalent to $cd(ca\gamma) < c(dd\delta) =
c\delta$.  It is clear that $cd(ca\gamma) < c\delta$ implies
$ca\gamma < c\delta$, because $ca\gamma \le cd(ca\gamma)$.  On the
other hand, if $ca\gamma < c\delta$, then $ca\gamma < \crit{cd}$,
so $cd(ca\gamma) = ca\gamma < c\delta$.  \QED

\proclaim{Lemma 6.9} If $a \tLequiv\gamma/b$ and\/ $[a,\delta] \le^*
H(\gamma)$, then\/ $[a,\delta] = [b,\delta]$.  \endproclaim

\demo{Proof} It is enough to show that $[b,\delta] \le^* [a,\delta]$,
since then one can interchange $a$ and $b$.  Fix $d$ such that
$\crit d = \delta$.  By the preceding lemma, we have $da\delta < d\gamma$.
This allows us to conclude from $da \tLequiv d\gamma/ db$ that
$da\delta = db\delta$; since Lemma~6.4(a) gives
$dd(db\delta) < d(db\delta)$, we get $dd(db\delta) < d(da\delta)$,
so $[b,\delta] \le^* [a,\delta]$, as desired.  \QED

We are now ready to prove the remaining property of $\tLequiv*/$: if $a
\tLequiv*\,\gamma^*/ b$ and $a[c,\rho] <^* \gamma^*$, then $a[c,\rho] =
b[c,\rho]$.  Fix $\delta$ such that $\gamma^* \le^* H(\delta)$ and $a
\tLequiv\delta/ b$.  We have $[ac,a\rho] <^* H(\delta)$, so the statement
preceding Lemma~6.8 gives $H(a\rho) < H(\delta)$.  Since $H$ is
order-preserving, we have $a\rho < \delta$.  Therefore, $a\rho = b\rho$,
so, using $ac \tLequiv\delta/ bc$ and Lemma~6.9, we get $a[c,\rho] =
[ac,a\rho] = [bc,a\rho] = [bc,b\rho] = b[c,\rho]$.

We have now completed the proof that $\E$ and $\O^*$, together with the
starred operations and relations, form a two-sorted embedding algebra.
Also, we have a canonical order-preserving map $H$ from $\O$ to
$\O^*$, and it is easy to check that $H$ sends all of the operations and
relations to their starred equivalents; hence, $(\E,\O)$ is isomorphic to a
subalgebra of $(\E,\O^*)$, so $(\E,\O^*)$ is isomorphic to an extension of
$(\E,\O)$.  It now remains to define the additional operations of an extended
two-sorted embedding algebra for $(\E,\O^*)$.

Since we want the pair $[b,\gamma]$ to represent $\appless b\gamma$, the
formula $\appless a{\appless b\gamma} = \appless{(a\circ b)}\gamma$
indicates that
we should define $\appless a{[b,\gamma]}$ to be $[a{\circ}b,\gamma]$.  The fact
that this is a valid definition (i.e., it does not depend on the choice
of a representative $(b,\gamma)$ for the equivalence class $[b,\gamma]$)
follows from the next lemma.

\proclaim{Lemma 6.10} If\/ $[b,\gamma] \le^* [b',\gamma']$, then\/
$[a{\circ}b,\gamma] \le^* [a{\circ}b',\gamma']$. \endproclaim

\demo{Proof} Fix $c$ and $c'$ such that $\crit c = \gamma$ and $\crit{c'} =
\gamma'$.  Since $[b,\gamma] \le^* [b',\gamma']$, we have $cc'(cb\gamma) \le
c(c'b'\gamma')$; applying $c(c'a)$ to this gives $c(c'a)(cc'(cb\gamma))
\le c(c'a)(c(c'b'\gamma'))$.  But $$c(c'a)(cc'(cb\gamma)) =
cc'(ca)(cc'(cb\gamma)) = cc'(ca(cb\gamma)) = cc'(c(a\circ b)\gamma)$$ and
$c(c'a)(c(c'b'\gamma')) = c(c'a(c'b'\gamma')) = c(c'(a\circ
b')\gamma')$, so we have $cc'(c(a\circ b)\gamma) \le c(c'(a\circ
b')\gamma')$ and hence $[a{\circ}b,\gamma] \le^* [a{\circ}b',\gamma']$. \QED

So $\appless a{[b,\gamma]}$ is well-defined.  The next lemma shows that this
definition
matches the original motivation.

\proclaim{Lemma 6.11} For all $a$ and\/ $\gamma$, $\appless a{H(\gamma)}=
[a,\gamma]$.  \endproclaim

\demo{Proof} Fix $c$ such that $\crit c = \gamma$; then $\appless a{H(\gamma)}=
[a{\circ}c,\gamma]$.  We have $\crit{ac} = a\gamma$, so $ac \tLequiv
a\gamma/ \id$, so $a \tLequiv a\gamma/ ac \circ a = a\circ c$.  From
$\gamma < c\gamma$, we get $ca\gamma < ca(c\gamma) =
c(a\gamma)$, so Lemma~6.8 gives $[a,\gamma] \le^* H(a\gamma)$.  Therefore,
Lemma~6.9 gives $[a,\gamma] = [a{\circ}c,\gamma]$, as desired. \QED

We now verify that this definition of $\appless a{\gamma^*}$ satisfies the
first five axioms listed in Definition~6.1.  Let $\gamma^* = [c,\rho]$
and $\delta^* = [d,\sigma]$.  The first two axioms are proved by
simple computations:
$$\gather
a(\appless b{\gamma^*}) = a[b{\circ}c,\rho] = [ab{\circ}ac,a\rho] =
\appless{ab}{[ac,a\rho]} = \appless{ab}{a\gamma^*}, \\
\appless a{\appless b{\gamma^*}} = \appless a{[b{\circ}c,\rho]} =
[a{\circ}b{\circ}c,\rho] =
\appless {(a\circ b)}{\gamma^*}.
\endgather$$

The next two axioms are equivalent to: $\gamma^* \le^* \delta^*$ if and
only if $\appless a{\gamma^*}\le^* a\delta^*$.  To prove this, fix $r$ and
$s$ such that $\crit r = \rho$ and $\crit s = \sigma$; then
$$\align
\appless a{\gamma^*}\le^* a\delta^*
&\iff [a{\circ}c,\rho] \le^* [ad,a\sigma] \\
&\iff r(as)(r(a\circ c)\rho) \le r(as(ad)(a\sigma) \\
&\iff ra(rs)(ra(rc\rho)) \le r(a(sd\sigma)) \\
&\iff ra(rs(rc\rho)) \le ra(r(sd\sigma)) \\
&\iff rs(rc\rho) \le r(sd\sigma) \\
&\iff \gamma^* \le^* \delta^*.
\endalign$$

For the fifth axiom, suppose $a \tLequiv*\,\gamma^*/ b$ and
$\appless a{\delta^*}\le^* \gamma^*$.  Find $\eta$ such that
$a \tLequiv\eta/ b$
and $\gamma^* \le^* H(\eta)$; then $\appless a{\delta^*}\le^* H(\eta)$
and $a\circ d \tLequiv\eta/ b\circ d$, so Lemma~6.9 gives
$\appless a{\delta^*}= [a{\circ}d,\sigma] = [b{\circ}d,\sigma] =
\appless b{\delta^*}$.

It remains to find a suitable definition for the cofinality function.
Since $[a,\gamma]$ is supposed to represent $\appless a\gamma$, where $\gamma$
is a critical point and hence regular, we define $\cf{}[a,\gamma]$ to
be $H(\gamma)$.  As usual, we need a lemma showing that this does
not depend on the choice of a representative for the equivalence class
$[a,\gamma]$.

\proclaim{Lemma 6.12} If\/ $\gamma \ne \delta$, then\/ $[a,\gamma] \ne
[b,\delta]$.  \endproclaim

\demo{Proof} We may assume $\gamma < \delta$.  Fix $c$ and $d$ such that
$\crit c = \gamma$ and $\crit d = \delta$; then $\crit{dc} = d\gamma = \gamma$
and $\crit{dcd} = dc\delta \ge \delta > \gamma$.  We can use $dc$ instead of
$c$ when comparing $[a,\gamma]$ with $[b,\delta]$:
$[a,\gamma] \le^* [b,\delta]$ iff $dcd(dca\gamma)
< dc(db\delta)$. Now the assumption
that $[a,\gamma] = [b,\delta]$ leads to a contradiction as follows:
$$\alignat 2
dcd(dca\gamma)
&< dc(db\delta) &&\qquad\qquad \text{since } [a,\gamma] \le^* [b,\delta] \\
&< d(ca\gamma) &&\qquad\qquad \text{since } [b,\delta] \le^* [a,\gamma] \\
&= d(ca)(d\gamma) &&\\
&= dc(da)\gamma &&\\
&= dcd(dca)(dcd\gamma) &&\\
&= dcd(dca\gamma). &&
\endalignat\nopagebreak$$ \QED

It is now trivial to verify the axioms $\cf(\critfunc^*(a)) = \critfunc^*(a)$,
$\cf(\appless a{\gamma^*}) = \cf\gamma^*$, and $\cf(a\gamma^*) =
a(\cf\gamma^*)$.

The last two axioms can actually be deduced from the other axioms when
the ordinal~$\gamma$ is of the form $\appless b\delta$ where $\delta$ is a
critical point; since every element of $\O^*$ has this form, this will
suffice here.  The law $\cf\gamma \le \gamma$ follows from
Proposition~6.2\therosteritem1, since $\cf\gamma = \cf\delta = \delta$.
Now suppose $a$ does not move $\delta = \cf\gamma$; then
$\appless a{\appless b\delta}
= \appless{(a\circ b)}\delta = \appless{(ab \circ a)}\delta =
\appless{ab}{\appless a\delta}$
and $a(\appless b\delta) = \appless{ab}{a\delta}$, and these two ordinals
are equal because $\delta \le \appless a\delta \le a\delta = \delta$.

This completes the proof of Theorem~6.3. \QED

Theorem~6.3 can be used to transfer various arguments from the context
of elementary embeddings to that of two-sorted embedding algebras.
One example is the following result, which Laver proved for
elementary embeddings (Theorem 2.14).

In a two-sorted embedding algebra, let $j \neq \id$ be some 
embedding, and let $A_j$ be the set of embeddings generated from
$j$ by the operation $\cdot$ (so each $a \in A_j$ is given by
a word in $W_\cala$).

\proclaim{Theorem 6.13} Assume that the set of all critical points
of elements of $A_j$ has order type $\omega$.
If $a$ and $b$ are distinct elements of $A_j$, then there is
a critical point\/ $\gamma$ such that $a(\gamma) \ne b(\gamma)$. \endproclaim

\demo{Proof} We may assume that
all ordinals in the algebra are critical points; otherwise, just move to
the subalgebra comprising all embeddings and all critical points.
Apply Theorem~6.3 to construct an extended two-sorted embedding algebra
which is an extension of the given algebra.
We now follow the proof of
Theorem~13 from Laver~\cite{\LavOAEER}; every step except one in this proof
uses only properties of the extended ordinals which are listed in~6.1
and~6.2, and hence works in the same way here.  The one exception is 
the use of the fact that a certain increasing sequence of critical points
is cofinal in the set of all critical points of $A_j$; we have made this
fact an assumption of the theorem.
The result is that, in the extended algebra, there
exists a critical point $\gamma$ such that $a\gamma \ne b\gamma$.
But all critical points in the extended algebra are critical points
in the original algebra (since the same holds for embeddings), so
so we have the desired result in the original algebra. \QED

\bigpagebreak
\head 7. Construction of an embedding algebra \endhead

In this section, we will prove one direction of Theorem~5.2 by showing
how to construct an embedding algebra under the assumption that
$A_\infty$ is free (and hence all of the statements in Theorem~4.4
hold).

We will first construct a two-sorted embedding algebra.
The embedding set $\E$ will be $P_\infty \cup \{0\}$, while the
ordinal set $\O$ will be $\natnum$.  The operations $\cdot$ and $\circ$
on $\E$ will of course be those obtained from $P_\infty$, and $0$ will
be the identity in $\E$.

We note that 4.4(vi) implies the stronger statement that, for every
$a \in W_\calp$, there is an $n$ such that $[a]_n \ne 0$.  To see
this, use Lemma~2.11 to find a word in $W_\calp$ of the form
$a_1 \circ \dots \circ a_k$ ($a_1,\dots,a_k \in W_\cala$) which is
equivalent to $a$.  By 4.4(vi), there exists $n$ so large that
$[a_i]_n \ne 0$ for all~$i$; then formula (3.3) implies that
$[a]_n \ne 0$.

For each $a \in W_\calp$, define the function $e_a \colon \natnum \to
\natnum$ as follows:  for each $n \in \natnum$, let $e_a(n)=s(a\cdot 2^n)$.
In other words, $e_a(n)$ is the
largest $m$ such that $[a \cdot 2^n]_m = 0$.  (By the strengthened
4.4(vi), there is a
largest such $m$ for each $n$.)
If $a = b$ in $P_\infty$, then $[a]_m =
[b]_m$ and $[a \cdot 2^n]_m = [a]_m * [2^n]_m = [b]_m * [2^n]_m =
[b\cdot 2^n]_m$ for all $n$ and $m$; hence, $e_a = e_b$.  It therefore
makes sense to write $e_a$ for $a \in P_\infty$. This will give the
desired application function from $\E\times\O$ to $\O$, so we will
sometimes write $a(n)$ for $e_a(n)$ (but not $an$, as this might be
confused with $a\cdot n$). Define $e_0$ to be the identity function on
$\natnum$.

For any $a \in P_\infty$, we can apply Proposition~3.2 to
show that, if $[a \cdot 2^n]_m = 0$, then $[a \cdot 2^{n+1}]_{m+1} =
0$; it follows that the function $e_a$ is strictly increasing.  (This is
obviously true for $e_0$ as well.)  Now induction gives $e_a(n) \ge n$
for all $n$.

Next, we prove that $e_{a \circ b} = e_a \circ e_b$ (i.e., the algebra
operation $\circ$ represents composition).  This follows from the
corollary to Lemma~4.3:
$$
e_{a\circ b}(n) = s((a\circ b)\cdot 2^n) = s(a\cdot(b\cdot 2^n)) =
s(a \cdot 2^{s(b \cdot 2^n)}) = s(a \cdot 2^{e_b(n)}) = e_a(e_b(n)).
$$

For $a \in P_\infty$, define $\crit a$ to be the largest $m$ such that
$[a]_m = 0$, as given by the strengthened~4.4(vi).
(We will see later that this is the
critical point of $e_a$.)  It follows that $[a]_{m+1} = 2^m$, so $[a \cdot
2^m]_{m+1} = [2^m\cdot 2^m]_{m+1} = 0$.  (For the last equality,
see Proposition~3.2.)  This proves that $a(\crit
a) > \crit a$.

We now define $\tLequiv N/$ for $N \in \O$ by: $a \tLequiv N/ b$ iff
$[a]_N = [b]_N$.  The fact that the algebra $P_N$ satisfies (\Lavlaws)
immediately implies most of the desired properties of $\tLequiv N/$.  In
particular, if $a \tLequiv N/ b$ and $a(m) < N$, then $[a \cdot 2^m]_N$
is nonzero, and $[b \cdot 2^m]_N$ must have the same nonzero value, so
we find that $a(m) = b(m)$.  The only remaining property that is
nontrivial is coherence, for which we argue as follows. Suppose $[a]_N =
[b]_N$ and $M = e_c(N)$; we must show that $[ca]_M = [cb]_M$.  The
definition of $M$ implies that $[c\cdot 2^N]_M = 0$, so the period
of $c$ in $P_M$ divides $2^N$.  But $[a]_N = [b]_N$, so $[a]_M$ and $[b]_M$
are congruent modulo $2^N$; therefore, $[ca]_M = [cb]_M$, as desired.

This completes the construction of the two-sorted embedding algebra.
The point of constructing this intermediate algebra is that it allows us
to apply Theorem~6.13 to conclude that, if $a \ne b$ in $A_\infty$,
then $e_a \ne e_b$.

We now construct an embedding algebra as follows.  Let $A = \{e_a
\colon a \in A_\infty\}$.
Define the operation~$\cdot$ on $A$ by the formula
$e_a \cdot e_b = e_{ab}$; this definition is valid because the mapping
from $a$ to $e_a$ is one-to-one.  It is clear 
that $A$ is generated from the single
function $e_1$ by the operation $\cdot$.

Proposition~5.4\therosteritem1 implies that the critical point of $e_a$
is equal to the number $\crit a$ defined above.
\comment
We next want to show that the relation $\tLequiv N/$ defined above
coincides with the relation~$\Lequiv N/$ defined for $A_j$ in
Definition~5.1; more precisely, $a \tLequiv N/ b$ if and only if $e_a
\Lequiv N/ e_b$.  To show that $a \tLequiv N/ b$ implies $e_a \Lequiv N/
e_b$, first note that, if $a \tLequiv N/ b$, then $ac_1\dotsm c_k
\tLequiv N/ bc_1\dotsm c_k$; therefore, it suffices to show that $a
\tLequiv N/ b$ implies $e_a \Lrestrict N = e_b \Lrestrict N$.  But this
is just a restatement of the fact that, if $a \tLequiv N/ b$ and $a(m) <
N$, then $a(m) = b(m)$.

Conversely, suppose that $a \not\tLequiv N/ b$, so $[a]_N \ne [b]_N$.  Let
$k$ be such that $[a]_N+k = 2^N$, and let $c_1 = \dots = c_k = 1$; then
$[ac_1\dotsm c_k]_N = 0 \ne [bc_1\dotsm c_k]_N$.  Hence, if $m
= \crit{bc_1\dotsm c_k}$, then $m < N$ and $e_ae_{c_1}\dotsm e_{c_k}(m)
= m < e_be_{c_1}\dotsm e_{c_k}(m)$.  Therefore, $e_ae_{c_1}\dotsm
e_{c_k}(m) \Lrestrict N \ne e_be_{c_1}\dotsm e_{c_k}(m) \Lrestrict N$,
so $e_a \not\Lequiv N/ e_b$.
\endcomment
Given this, it is easy to see that $A$ satisfies the
axioms of an an embedding algebra by using the corresponding properties
of the two-sorted embedding algebra.  This completes the construction.

\bigpagebreak
\head 8. Uniqueness of embedding algebras \endhead

In this section, we will 
prove the following uniqueness result for monogenic embedding algebras.

\proclaim{Theorem 8.1} {\rm(a)} If $(A,{\cdot})$ is a monogenic embedding
algebra for which every natural number is a critical point, then
$(A,{\cdot})$ is isomorphic to the embedding algebra constructed
from~$P_\infty$ in the preceding section.

{\rm(b)} If $(\E,\O;{\cdot},{\circ},\dotsc)$ is a monogenic two-sorted
embedding algebra in which the ordinals have order type~$\omega$ and
every ordinal is a critical point, then it is isomorphic to the
two-sorted embedding algebra constructed from~$P_\infty$ in the
preceding section. \endproclaim

Along the way, we will show that, if a nontrivial
embedding algebra (or a nontrivial two-sorted embedding algebra
with ordinals of order type~$\omega$) exists, then
4.4(vi) holds, and hence $A_\infty$~is free, thus
completing the proof of Theorem~5.2.
Most of the arguments in this
section are adapted from Laver~\cite{\LavOAEER}.

If an embedding algebra satisfies the hypotheses of
Theorem~8.1(a), then, as noted after the proof of Proposition~5.5, we
can expand/extend it to a two-sorted
embedding algebra as hypothesized in Theorem~8.1(b).
So let us assume we have such a two-sorted embedding algebra.
Let $j$~be the generating embedding, and let $A_j$ be the set of
embeddings generated from~$j$ using $\cdot$ alone.
As noted in section~5, the set of non-identity embeddings is closed
under~$\cdot$, so every element of~$A_j$ has a critical point.
For any $a \in W_\cala$, let $j_a$ be the result of replacing each $1$ in the
expression~$a$ with~$j$.  (Note that $A_j = \{j_a \colon a \in W_\cala\}$.)  In
particular, since we identified positive integers
with words in $W_\cala$, we have an embedding $j_m$ for each $m>0$, and
$j_1=j$;
also, we let $j_0=\id$.

Let $\gamma_n$ be the critical point of $j_{2^n}$.
Recall that, for any $a \in W_\cala$, $[a]_n$ is defined to be the
result of evaluating $a$ in $A_n=\{0,1,\dots,2^n-1\}$.

\proclaim{Proposition 8.2} For any $a \in W^*_\cala$, $j_a \tLequiv\gamma_n/
j_{[a]_n}$. \endproclaim

\demo{Proof}  Since $j_{2^n} \tLequiv\gamma_n/ \id = j_0$, it does
not matter whether we work with $A_n$ or $A'_n=\{1,\dots,2^n\}$.
Clearly the proposition holds for $a = 0$.  We will show that, for any $b,c
\in A'_n$, $j_bj_c \tLequiv\gamma_n/ j_{b *_n c}$; given this, an easy
induction on $a \in W_\cala$ yields the proposition.

The proof of $j_bj_c \tLequiv\gamma_n/
j_{b *_n c}$ is by double induction, downward on $b$ and upward on
$c$.  For $b=2^n$, we have $j_{2^n}j_c \tLequiv\gamma_n/ \id \cdot j_c =
j_c = j_{2^n *_n c}$.  The case $b < 2^n$, $c=1$ is also trivial: $j_bj_1
= j_{b+1} = j_{b *_n 1}$.  Finally, for $b,c < 2^n$,
$$j_bj_{c+1} = j_b(j_c j) = (j_bj_c)(j_b j) \tLequiv\gamma_n/
j_{b *_n c}j_{b+1} \tLequiv\gamma_n/
j_{(b *_n c) *_n (b+1)} = j_{b *_n (c+1)}.$$
(The induction hypothesis can be used in the second-to-last step because
$b *_n c > b$.)  This completes the induction. \QED

\proclaim{Proposition 8.3} For all $n$, $\gamma_n <
\gamma_{n+1}$; also, for all $m > 0$, $\crit{j_m} = \gamma_k$ where\/ $2^k$ is
the largest power of\/ $2$ dividing $m$. \endproclaim

\demo{Proof} By induction on $N$, we show that these statements are true
for $n < N$ and $m < 2^N$.  The case $N = 0$ is vacuous.  Suppose now that
the assertion is true for $N$; we will prove it for $N+1$.  We know that
$\gamma_0 < \gamma_1 < \dots < \gamma_N$.  By definition, $\crit{j_m} =
\gamma_N$ if $m = 2^N$.  If $2^N < m < 2^{N+1}$, then Proposition 8.2
implies that $j_m \tLequiv\gamma_N/ j_{m-2^N}$,
so, if $2^k$ is the largest power of $2$ dividing $m-2^N$, then $2^k$ is also
the largest power of $2$ dividing $m$, and $\crit{j_{m-2^N}} = \gamma_k
< \gamma_N$, so $\crit{j_m} = \gamma_k$.  This means that the embeddings
$j_m$ for $2^N < m < 2^{N+1}$ all have critical points below $\gamma_N$,
and hence, by Proposition~5.7,
$j_m(\gamma_N) > \gamma_N$; let $\theta > \gamma_N$ be the least
of these values $j_m(\gamma_N)$.  Now coherence gives
$j_{m+1} = j_m j \tLequiv\theta/ j_m(j_{2^N} j)$ for all~$m$ in this range, so
$$\align j_{2^{N+1}} &\tLequiv\theta/ j_{2^{N+1}-1} (j_{2^N} j) \\
&\tLequiv\theta/ j_{2^{N+1}-2} (j_{2^N} j) (j_{2^N} j) \\
&\tLequiv\theta/ \dotso \\
&\tLequiv\theta/ j_{2^N+1} (j_{2^N} j) \dotsm (j_{2^N} j) =
j_{2^N}j_{2^N}. \endalign$$
Since $\crit{j_{2^N}j_{2^N}} = j_{2^N}(\gamma_N) > \gamma_N$ and
$\gamma_N < \theta$, we must have $\crit{j_{2^{N+1}}} > \gamma_N$.  This
completes the induction. \QED

We can now show that 4.4(vi) holds, and hence $A_\infty$ is free,
as follows: Suppose
$a \in W_\cala$.  Since the sequence of critical points $\gamma_n$ is
strictly increasing, and the ordinals have order type~$\omega$,
there must be an $n$ such that $\crit{j_a} < \gamma_n$.
Then $j_a \not\tLequiv\gamma_n/ \id = j_0$, so
we must have $[a]_n \ne 0$.  This completes the proof of Theorem~5.2.

\medpagebreak

Proposition 8.3 implies that $j_m \not\tLequiv\gamma_n/ \id$ for $1 \le m <
2^n$ (because $\crit{j_m} < \gamma_n$).  
Consequently, we have $j_m \not\tLequiv\gamma_n/ j_{m'}$ for
$1 \le m < m' \le 2^n$; if this were not so, then one could apply $j_m$ and
$j_{m'}$ to $j$ $2^n-m'$ times to get $j_{m+2^n-m'} \tLequiv\gamma_n/ j_{2^n}
\tLequiv\gamma_n/ \id$,  a contradiction.

It follows that the
mapping $j_a/{\tLequiv\gamma_n/} \mapsto [a]_n$ from
$A_j/{\tLequiv\gamma_n/}$ to~$A_n$ 
is bijective and preserves the operation $\cdot$, so it is an isomorphism.
These mappings commute with the canonical projections from
$A_j/{\tLequiv\gamma_{n+1}/}$ to $A_j/{\tLequiv\gamma_n/}$ and from
$A_{n+1}$ to~$A_n$, so they give a mapping from $A_j$ to the inverse
limit of the algebras $A_n$; clearly this mapping sends the generator of
$A_j$ to the generator of $A_\infty$, so we have a mapping $f$ from
$A_j$ onto $A_\infty$.  Since $f$ preserves $\cdot$, and since
$A_\infty$ is free, $f$ must be an isomorphism between $A_j$ and
$A_\infty$.

For any $a \in W_\cala$, if $n$~is so large that
$\crit{j_a} < \gamma_n$, then Proposition~8.2 gives
$j_a \tLequiv\gamma_n/ j_{[a]_n}$, and $[a]_n$~must be nonzero,
so, by Proposition~8.3, $\crit{j_a} = \gamma_k$ where
$2^k$ is the largest power of~$2$ dividing~$[a]_n$.  This~$k$
is just~$s(a)$.  Also, for any~$m$, we get
$$j_a(\gamma_m) = j_a(\crit{j_{2^m}}) = \crit{j_aj_{2^m}}
= \crit{j_{a\cdot 2^m}} = \gamma_{s(a \cdot 2^m)} = \gamma_{e_a(m)},$$
where $e_a$~is as defined in section~7.
Finally, for any $a,b \in W^*_\cala$, we have
$$j_a \tLequiv\gamma_n/ j_b \iff
j_{[a]_n} \tLequiv\gamma_n/ j_{[b]_n} \iff [a]_b = [b]_n.$$

Therefore, the structure of $A_j$ is determined completely
except for the possible existence of ordinals which are not critical points.
(Even for these, the equivalence relation $\tLequiv \delta/$ is determined; the
argument of the preceding paragraph shows that, if $\gamma_{n-1} < \delta \le
\gamma_n$, then $j_a \tLequiv \delta/ j_b$ if and only if $[a]_n = [b]_n$.)
In the situation of Theorem~8.1(a), there are no such extra ordinals,
and we have $\gamma_n = n$ for all~$n$; we can now see that the
structure of~$A_j$ (which is just a copy of the original embedding
algebra~$A$) exactly matches the structure defined in section~7
from~$P_\infty$.  So Theorem~8.1(a) is proved.

Now, in the situation of Theorem~8.1(b), let $P_j$~be the set of embeddings
generated from~$j$ using both $\cdot$ and~$\circ$.  (Since the
algebra is generated by~$j$, this is all embeddings except~$\id$.)
In order to show that composition here matches the structure from
section~7, we use the following result.

\proclaim{Proposition 8.4} If $a,b \in A_j$ and $n \in \natnum$,
then there is $c \in A_j$ such that $a \circ b \tLequiv\gamma_n/ c$.
\endproclaim

\demo{Proof} We may assume that $\crit a > \crit b$; otherwise, replace
$a$ and $b$ with $ab$ and $a$ (using $a\circ b = ab \circ a$).  Now let
$a_0 = b$, $a_1 = a$, and $a_i = a_{i-1}a_{i-2}$ for $i \ge 2$.  Induction
gives $a_{i+1} \circ a_i = a \circ b$,
$\crit b = \crit{a_0} = \crit{a_2} = \crit{a_4} = \dotso$, and
$\crit a = \crit{a_1} < \crit{a_3} < \crit{a_5} < \dotso$.  Since the
sequence $\crit{a_{2i+1}}$ is a strictly increasing sequence
of critical points, and the set of all critical points has order type~$\omega$,
there must be an odd $i$ such that $\crit{a_i} \ge
\gamma_n$; this gives $a \circ b = a_i \circ a_{i-1} \tLequiv\gamma_n/
a_{i-1}$, so we can let $c = a_{i-1}$. \QED

It follows that, if $a,b \in A_n$, then there is $c \in A_n$ such that
$j_a \circ j_b \tLequiv\gamma_n/ j_c$; we know from the above results that
this $c$ is unique.  To determine what $c$ is, note that $(j_a \circ j_b)j
\tLequiv\gamma_n/ j_cj$, so $j_{a *_n (b+1)} = j_{c+1}$, so $a *_n (b+1)
= c+1$, where the additions are performed modulo~$2^n$ in~$A_n$; hence,
$c = (a *_n (b+1)) - 1 = a \circ_n b$.  We therefore have $j_a \circ j_b
\tLequiv\gamma_n/ j_{a \circ_n b}$ for $a,b \in A_n$; now, if we define
$j_a$ for $a \in W_\calp$ as we did for $a \in W_\cala$, then induction on
$a$ gives $j_a \tLequiv\gamma_n/ j_{[a]_n}$ for all
$a \in W_\calp$.  We can now argue as before that $P_j/{\tLequiv\gamma_n/}$
is isomorphic to $P_n$ and $P_j$ is isomorphic to $P_\infty$, so the
structure of~$P_j$ is unique except for the possible existence of ordinals
which are not critical points, and matches that from section~7.
This completes the proof of Theorem~8.1.

\medpagebreak

One can in fact construct an embedding algebra with numbers that are not
critical points, either by just duplicating every critical point or, less
trivially, by constructing the extended algebra in section 6 and then using
the method of section 7 to convert this to an embedding algebra.
(One can then modify the algebra further to get an
embedding algebra which does not satisfy $ab(a(n)) = a(b(n))$.)  For the
less trivial construction, one must observe that the ordinals in the
extended algebra
have order type $\omega$. To see this, note that if $a \tLequiv\gamma_n/ b$
and $\appless a\kappa < \gamma_n$, then $\appless a\kappa = \appless b\kappa$;
hence,
there are at most $n2^n$ extended ordinals below $\gamma_n$.

On the other hand, we now have a roundabout proof that, if there is
a nontrivial embedding algebra, then there is one in which
all natural numbers are critical points (and hence
$ab(a(n)) = a(b(n))$ holds), namely the one constructed from~$P_\infty$.
One would expect to be able to prove this directly, by simply
deleting the natural numbers which are not critical points
and relabeling the critical points as $0,1,2,\dots$.  However, it is
conceivable that distinct functions in the algebra are the same when
restricted to the critical points, so that $\cdot$~could fail to be
well-defined after the other numbers are deleted.  It turns out that
this does not happen in the monogenic case, but the authors
do not see a way to prove this without building up enough
structure to imitate Laver's proof of Theorem~2.14.

\bigpagebreak
\head 9. The strength of ``$A_\infty$ is free'' \endhead

As we recalled in section 1, Laver's proof of the irreflexivity of the
free left distributive algebra on one generator assumed
the existence of a nontrivial elementary embedding from~$V_\lambda$ to
itself; this is an extremely strong large cardinal hypothesis.
(Actually, Laver had noted that, since one only needs a bounded part of
$V_\lambda$ to talk about the finitely many embeddings mentioned while
comparing two given words in the free algebra, the assumption can be
reduced to the existence of an $n$-huge cardinal for each natural number~$n$.)
The possibility that the irreflexivity property was strong enough
to require large cardinal assumptions for its proof remained until
Dehornoy proved the property without such assumptions (in fact, using
only Primitive Recursive Arithmetic).

We now consider the statement ``$A_\infty$ is free'' and the equivalent
versions in Theorem~4.4.  These statements imply that $A_\infty$ is both
free and irreflexive, so the irreflexivity of the free algebra follows
immediately.  The purpose of this section is to show that the statement
``$A_\infty$ is free'' is strictly stronger than the statement ``the
free algebra is irreflexive,'' in the following sense:

\proclaim{Theorem 9.1} The statement ``$A_\infty$ is free'' is not
provable in Primitive Recursive Arithmetic. \endproclaim

Of course, we assume throughout that PRA is itself consistent.

\demo{Proof} It is a well-known result from proof theory (see
Sieg~\cite{\SieFA}) that the only recursive functions that can be proved to
be total using only PRA are the primitive recursive functions.
Therefore, to prove the theorem, it will suffice to show that
PRA${}+{}$4.4(vii) proves the totality of a recursive function $F$ which
is not primitive recursive.

For each natural number $n$, let $F(n)$ be the largest $m$ such that
$[u_n]_m = 0$, where $u_n$ is the word
$1\cdot(1\cdot(\dots(1\cdot1)\dots))$ with $n+1$ $1$'s.  It
follows from 4.4(vii) that $F$ is a total recursive function.  If the
functions $e_a$ are defined as in section~7, thus giving an embedding
algebra, then $F(n) = e_1^n(0)$, so $F$ is the critical sequence of the
mapping $e_1$. Since all natural numbers are critical points in this
embedding algebra, one can state that $F(n)$ is the number of critical
points below $e_1^n(0)$.

We now use the methods of Dougherty~\cite{\DouCPAEE}
for producing many critical
points.  That paper is written in terms of elementary embeddings, but it
is not hard to check that the only properties used in section~2 of that
paper are that each embedding gives a strictly increasing monotone
function on the ordinals and that, if $a$ and $b$ are two such
embeddings, then $\crit{ab} = a(\crit b)$ and $a(b\gamma) = ab(a\gamma)$
for all ordinals $\gamma$.  Hence, the main theorem of that paper, that
the number of critical points below $\kappa_n = j^n(\kappa_0)$
grows so rapidly with~$n$ that it cannot be primitive recursive, applies
to any nontrivial embedding algebra or two-sorted embedding algebra. (The
results in later sections of that paper use only the properties of an
extended two-sorted embedding algebra, so the stronger lower bounds
obtained there also apply to any nontrivial embedding algebra or two-sorted
embedding algebra.)  But, in the embedding algebra from section~7, the
number of critical points below~$\kappa_n$ is just $F(n)$ as defined
above, so $F$ is not primitive recursive. \QED

On the other hand, the freeness of $A_\infty$ follows from the existence
of a nontrivial elementary embedding $j \colon V_\lambda \to
V_\lambda$.  The proof of this (due to Laver) uses Theorem~2.13.
Given this theorem, we can
apply the arguments in section~8 to the monogenic two-sorted embedding
algebra obtained from $P_j$ to conclude that $A_\infty$ is free.
Laver (personal communication) has recently noted, and the authors
have confirmed, that one can use the method of proof of Theorem~2.13
while working with only an $n$-huge embedding, to get a correspondingly
weaker result; hence, the freeness of $A_\infty$ follows from the
existence of an $n$-huge cardinal for each natural number $n$.
(There is a level-by-level form of this result: if a $k$-huge
cardinal exists, then there is a natural number $n$ such that
$[u_k]_n \ne 0$.)

The proof of Theorem~9.1 showed that the assumption that $A_\infty$ is
free can be used to construct a particular function $F$ which grows too
rapidly to be primitive recursive.  It turns out that one cannot produce
any function growing much faster than $F$ from this assumption.  This
can be stated precisely as follows.

\proclaim{Proposition 9.2} Any recursive function which is provably
total in\/ {\rm PRA${}+{}$}``$A_\infty$ is free'' must grow more slowly
than $F_m$ for some $m$, where $F_0 = F$ and $F_{m+1}$ is the iteration
of~$F_m$ (starting at\/~$1$, say; that is, $F_{m+1}(n) = F_m^n(1)$).
\endproclaim

\demo{Proof} The proofs in Sieg~\cite{\SieFA} can be modified to give the
following extended version of the proof-theoretic result used earlier:

{\narrower\narrower\noindent If $P(n,m)$ is a primitive recursive predicate,
$f(n)$ is the least $m$ such that $P(n,m)$ holds, and $g$ is a recursive
function which is provably total in PRA${}+ \forall n \exists m P(n,m)$,
then $g$ can be obtained from $f$ and trivial functions (constants,
projections, and successor) by composition and primitive recursion.

}
The function $F$ can be used as $f$, since $P(n,m)$ can be defined to be
``$u_n \ne 0$ in $A_{m+1}$.''  Also, 4.4(vii) is a consequence of
PRA${}+ \forall n \exists m P(n,m)$. Therefore, any
recursive function~$g$ provably total from PRA${}+{}$``$A_\infty$ is
free'' must be obtainable from $F$ and trivial functions by the
operations of composition and primitive recursion.  Now the standard
proof by induction on the number of such operations used shows that
$g$ is below $F_n$ for some $n$. \QED

As a particular case of this, recall that there is a primitive recursive
algorithm for comparing two expressions $a$ and $b$, i.e., transforming
them into equivalent expressions $a'$ and $b'$ such that either $a' =
b'$ or one of $a',b'$ is a left subterm of the other.  Starting with
this, one can go through Laver's proof that any two distinct embeddings
must differ at a critical point, and verify that all of the steps are
primitive recursive.  Hence, assuming $A_\infty$~is
free, if $a$ and $b$ are members
of~$W_\cala$ such that $a \not\equiv_{\cala} b$,
and $n$ is least such that $j_a$ and $j_b$ differ at
critical point number $n$, then $n$ can be obtained from
$a$ and~$b$ by a function whose growth rate is comparable to that of
$F$.

\bigpagebreak
\head 10. Open problems and acknowledgments \endhead

There remain a number of open problems related to these algebras.
The main one, of course, is the exact strength of the statement
``$A_\infty$ is free''; the gap between ``more than PRA'' and
``there is an $n$-huge cardinal for each $n$''
is rather large.
One can also ask whether ``there is a nontrivial two-sorted
embedding algebra'' is as strong as ``there is a nontrivial
embedding algebra.''

It is still open whether Laver's result on distinguishing elementary
embeddings by their behavior on critical points (Theorem~2.14) can
be extended to~$P_j$.  If it can, by methods formalizable in an
extended two-sorted embedding algebra, then one can define a version
of embedding algebra which includes a composition operation, and the
existence of a nontrivial such algebra will still be equivalent to
``$A_\infty$ is free.''

Another area of interest is further extensions of the results in
section~6 to include more of the ordinals that can be defined
from elementary embeddings.  (Eventually one might hope to start
with the embedding algebra obtained from $A_\infty$ and construct
a larger structure including all of the important features of
the algebra obtained from an elementary embedding from $V_\lambda$
to itself.)  A natural next step is to try to define ordinals
of the form ``the least $\alpha$ such that $a(\alpha) \ge \gamma$''
for a given embedding $a$ and ordinal $\gamma$.  Such ordinals
seem to be closely tied to the inequality $aa(\gamma) \le a(\gamma)$:
the existence of the ordinals allows one to prove that the inequality
holds, and the authors can show under the assumption of the inequality
that there is a natural extension of a given monogenic two-sorted
embedding algebra in which all ordinals are critical points
to an algebra including such ordinals.  The authors do not
yet have a large-cardinal-free proof that the inequality holds
in the embedding algebra constructed from $A_\infty$, even assuming
that $A_\infty$ is free.

The authors would like to thank 
P. Dehornoy for discussing his work and suggesting further questions, 
R. Laver for showing us his unpublished results,
M. Rathjen for consultations about proof theory,
and J. Zapletal for pointing out the construction used in Proposition 2.12.

\enddocument